\documentclass{gtart}
\usepackage{graphicx, amsmath, amssymb, labelfig}
\nocolon

\def\ifplaintex{\expandafter\ifx\csname documentclass\endcsname\relax}

\ifplaintex 
\else
\fi

\expandafter\ifx\csname beginpicture\endcsname\relax
\expandafter\ifx\csname documentclass\endcsname\relax
\input pictex \else
\input prepictex \input pictex \input postpictex \fi\fi

\def\gt{{\mathsurround=0pt\it $\cal G\mskip-2mu$eometry \&\ 
$\cal T\!\!$opology}}        

\def\gtp{{\mathsurround=0pt\it $\cal G\mskip-2mu$eometry \&\ 
$\cal T\!\!$opology $\cal P\!$ublications}}  

\def\volumenumber#1{\def\thevolumenumber{#1}}
\def\papernumber#1{\def\thepapernumber{#1}}
\def\volumeyear#1{\def\thevolumeyear{#1}}

\def\pagenumbers#1#2{\def\startpage{#1}\def\finishpage{#2}}
\def\published#1{\def\publishdate{#1}}
\def\proposed#1{\def\theproposer{#1}}
\def\seconded#1{\def\theseconders{#1}}
\def\received#1{\def\receiveddate{#1}}

\def\accepted#1{\def\accepteddate{#1}}

\def\asciiaddress#1{\def\theasciiaddress{#1}}
\def\asciiemail#1{\def\theasciiemail{#1}}
\long\def\asciiabstract#1{\long\def\theasciiabstract{#1}}

\def\shortauthors#1{\def\theshortauthors{#1}}

\let\\\par\let\thevolumenumber\relax\let\thepapernumber\relax
\let\thevolumeyear\relax\let\thesamplenumber\relax\let\startpage\relax
\let\finishpage\relax\let\publishdate\relax\let\receiveddate\relax
\let\reviseddate\relax\let\accepteddate\relax\let\theasciititle\relax
\let\theasciiauthors\relax\let\theasciiaddress\relax
\let\theasciiabstract\relax
\let\theasciiemail\relax\let\theshortauthors\relax\let\theshorttitle\relax

\long\def\maketitlep{   

\count0=\startpage

\gt\hfill      
\beginpicture
\setcoordinatesystem units <0.33truein, 0.33truein> point at 2.2 0.9
\setplotsymbol ({$\cal G$})
\plotsymbolspacing=9truept
\circulararc 315 degrees from 0 1 center at 0 0
\setplotsymbol ({$\cal T$})
\circulararc 315 degrees from 1 -1 center at 1 0
\endpicture
\break
{\small\ifx\thesamplenumber\relax 
Volume \else Sample
\fi\thevolumenumber\ (\thevolumeyear)
\startpage--\finishpage\nl
Published: \publishdate}
\vglue 0.5truein plus 0.4fil minus 0.1truein

{\parskip=0pt\leftskip 0pt plus 1fil\def\\{\par\smallskip}{\ifplaintex\large
\else\Large\fi\bf\thetitle}\par\medskip}   

\vglue 0pt plus 0.1fil 

{\parskip=0pt\leftskip 0pt plus 1fil\def\\{\par}{\sc\theauthors}
\par\medskip}

\vglue 0pt plus 0.1fil 

{\small\parskip=0pt\let\newline\\
{\leftskip 0pt plus 1fil\def\\{\par}{\sl\theaddress}\par}
\expandafter\ifx\theemail\relax    
\relax\else\vglue 5pt plus 0.02fil minus 2pt\def\\{\stdspace{\rm 
and}\stdspace} 
\cl{Email:\stdspace\tt\theemail}\fi
\ifx\theurl\relax                  
\relax\else\vglue 5pt plus 0.02fil minus 2pt\def\\{\stdspace{\rm 
and}\stdspace}
\cl{URL:\stdspace\tt\theurl}\fi\par}

\vglue 7pt plus 0.3fil minus 3pt

{\bf Abstract}
\vglue 5pt plus 0.1fil minus 2pt

\theabstract

\vglue 7pt plus 0.3fil minus 3pt

{\bf AMS Classification numbers}\quad Primary:\quad \theprimaryclass

Secondary:\quad \thesecondaryclass

\vglue 5pt plus 0.3fil minus 2pt

{\bf Keywords:}\quad \thekeywords

\vglue 10pt plus 0.5fil minus 5pt

{\small  Proposed: \theproposer\hfill Received: \receiveddate\nl
Seconded: \theseconders\hfill 
\ifx\reviseddate\relax                         
Accepted: \accepteddate                        
\else
Revised: \reviseddate                          
\fi}
\eject
}       

\let\maketitlepage\maketitlep
\let\maketitle\maketitlepage

\font\phead=cmsl9 scaled 950
\font=cmsl9 scaled 1050
\font\pnum=cmbx10 scaled 913
\font\lnum=cmbx10 
\font\pfoot=cmsl9 scaled 950
\font=cmsl9 scaled 1050
\ifplaintex
\headline{\vbox to 0pt{\vskip -4.5mm\line{\small\phead\ifnum
\count0=\startpage ISSN 1364-0380 (on line)
1465-3060 (printed) \hfill {\pnum\folio}\else\ifodd\count0\def\\{ }%
\ifx\theshorttitle\relax\thetitle\else\theshorttitle\fi\hfill{\pnum\folio}
\else\def\\{ and }{\pnum\folio}\hfill\ifx\theshortauthors\relax\theauthors
\else\theshortauthors\fi\fi\fi}\vss}}
\footline{\vbox to 0pt{\vglue 0mm\line{\small\pfoot\ifnum\count0=\startpage
\copyright\ \gtp\hfill\else
\gt, Volume \thevolumenumber\ (\thevolumeyear)\hfill\fi}\vss
}}
\else
\makeatletter
\def\@oddhead{{\small\lhead\ifnum\count0=\startpage ISSN 1364-0380 (on line)
1465-3060 (printed) \hfill {\lnum\number\count0}\else\ifodd\count0
\def\\{ }\ifx\theshorttitle\relax \thetitle \else\theshorttitle\fi\hfill
{\lnum\number\count0}\else\def\\{ and }{\lnum\number\count0}
\hfill\ifx\theshortauthors\relax 
\theauthors\else\theshortauthors\fi\fi\fi}}\def\@evenhead{@oddhead}
\def\@oddfoot{\small\lfoot\ifnum\count0=\startpage\copyright\ \gtp\hfill\else
\gt, Volume \thevolumenumber\ (\thevolumeyear)\hfill\fi}
\def\@evenfoot{@oddfoot}
\makeatother
\fi

\newwrite\gtoutfile
\long\gdef\makeheadfile{  
{\def\\{, }
\immediate\openout\gtoutfile head.xxx
\immediate\write\gtoutfile{To: math@arxiv.org}
\immediate\write\gtoutfile{Subject: put}
\immediate\write\gtoutfile{--text follows this line--}
\immediate\write\gtoutfile{Proxy-for: \ifx\theasciiauthors\relax
\theauthors\else\theasciiauthors\fi <\ifx\theasciiemail\relax\theemail\else\theasciiemail\fi>}
\immediate\write\gtoutfile{\noexpand\\}
\immediate\write\gtoutfile{Authors: \ifx\theasciiauthors\relax
\theauthors\else\theasciiauthors\fi}
{\def\\{ }\immediate\write\gtoutfile{Title: \ifx\theasciititle\relax
\thetitle\else\theasciititle\fi}}
\immediate\write\gtoutfile{Subj-class: GT}
\immediate\write\gtoutfile{MSC-class: \theprimaryclass\ifx\thesecondaryclass\relax\else, \thesecondaryclass\fi}
\immediate\write\gtoutfile{Journal-ref: Geom. Topol. \thevolumenumber
(\thevolumeyear) \startpage-\finishpage}
\immediate\write\gtoutfile{Comments: Published in Geometry and Topology at}
\immediate\write\gtoutfile{    http://www.maths.warwick.ac.uk/gt/GTVol\thevolumenumber/paper\thepapernumber.abs.html}
\immediate\write\gtoutfile{\noexpand\\}
\immediate\write\gtoutfile{}
\ifx\theasciiabstract\relax
\immediate\write\gtoutfile{\theabstract}\else
\immediate\write\gtoutfile{\theasciiabstract}\fi
\immediate\write\gtoutfile{}
\immediate\write\gtoutfile{\noexpand\\}
\immediate\write\gtoutfile{}
\immediate\write\gtoutfile{<uuencoded .tar.gz file here>}
\immediate\write\gtoutfile{}
\immediate\closeout\gtoutfile}}  

\def\maketitlepage{\maketitlep\makeheadfile}
\let\maketitle\maketitlepage

\volumenumber{4}\papernumber{8}\volumeyear{2000}
\pagenumbers{243}{275}
\proposed{Cameron Gordon}
\seconded{Robion Kirby, David Gabai}
\received{17 January 2000}
\accepted{18 September 2000}
\published{3 October 2000}

\newtheorem{theorem}{Theorem}[section]
\newtheorem{proposition}[theorem]{Proposition}
\newtheorem{lemma}[theorem]{Lemma}
\newtheorem{claim}[theorem]{Claim}
\newtheorem{corollary}[theorem]{Corollary}
\theoremstyle{definition}
\newtheorem{definition}[theorem]{Definition}

\newtheorem{remark}[theorem]{Remark}

\numberwithin{equation}{section}

\newcommand{\aaa}{\mbox{$\alpha$}}

\newcommand{\bbb}{\mbox{$\beta$}}

\newcommand{\gam}{\mbox{$\gamma$}}

\newcommand{\bdd}{\mbox{$\partial$}}

\newcommand{\Ddd}{\mbox{$\Delta$}}

\newcommand{\lam}{\mbox{$\lambda$}}

\newcommand{\inter}{\mbox{${\rm int}$}}

\begin{document}

\title{Levelling an unknotting tunnel}

\author{Hiroshi Goda\\Martin Scharlemann\\Abigail Thompson}
\shortauthors{Goda, Scharlemann and Thompson}

\address{Graduate School of Science and Technology, Kobe University\\
Rokko, Kobe 657-8501, Japan\\\smallskip\\
Mathematics Department, University of California\\
Santa Barbara, CA 93106, USA\\\smallskip\\
Mathematics Department, University of California\\
Davis, CA 95616, USA\\\smallskip\\
{\rm Email:}\qua\tt goda@math.kobe-u.ac.jp, mgscharl@math.ucsb.edu\\
thompson@math.ucdavis.edu}

\asciiaddress{Graduate School of Science and Technology, Kobe University\\
Rokko, Kobe 657-8501, Japan\\
Mathematics Department, University of California\\
Santa Barbara, CA 93106, USA\\
Mathematics Department, University of California\\
Davis, CA 95616, USA}

\asciiemail{goda@math.kobe-u.ac.jp, mgscharl@math.ucsb.edu,
thompson@math.ucdavis.edu}

\begin{abstract}
It is a consequence of theorems of Gordon--Reid \cite{GR} and Thompson
\cite{T} that a tunnel number one knot, if put in thin position, will
also be in bridge position.  We show that in such a thin
presentation, the tunnel can be made level so that it lies in a
level sphere.  This settles a question raised by Morimoto \cite{M1},
who showed that the (now known) classification of unknotting tunnels
for 2--bridge knots would follow quickly if it were known that any
unknotting tunnel can be made level.
\end{abstract}
\asciiabstract{It is a consequence of theorems of Gordon-Reid [Tangle
decompositions of tunnel number one knots and links, J. Knot Theory
and its Ramifications, 4 (1995) 389-409] and Thompson [Thin position
and bridge number for knots in the 3-sphere, Topology, 36 (1997)
505-507] that a tunnel number one knot, if put in thin position, will
also be in bridge position.  We show that in such a thin presentation,
the tunnel can be made level so that it lies in a level sphere.  This
settles a question raised by Morimoto [A note on unknotting tunnels
for 2-bridge knots, Bulletin of Faculty of Engineering Takushoku
University, 3 (1992) 219-225], who showed that the (now known)
classification of unknotting tunnels for 2-bridge knots would follow
quickly if it were known that any unknotting tunnel can be made
level.}



\maketitlepage

\section{Background}

In \cite{G}, Gabai introduced the notion of {\em thin position} for
a knot in $S^3$.   Choose a height function
$h\co S^3-\{x,y\}=S^2\times \mathbb R \rightarrow\mathbb R$, the
projection to the second factor.  Informally, a knot $K$ is in thin
position with respect to $h$ if the number of points in which it
intersects the level spheres of $h$ has been minimized.  More
formally, let
$S(t)=h^{-1}(t)$ and call it the sphere at {\it height} or {\it
level} $t$. Put $K$ in general position with respect to $h$ and choose
heights $k_1 < k_2 < \ldots < k_n$ between each successive pair of
critical heights of $K$.  Define the {\em width} $W(K)$ of this
imbedding of $K$ in $S^3$ to be the integer   $\Sigma_i|S(k_i) \cap
K|$.

\begin{definition} \label{def:knotthin} $K \subset S^3$ is in {\em
thin position} (with respect to the height function $h$) if $W(K)$
cannot be reduced by an isotopy of $K$ in $S^3$.
\end{definition}

\begin{definition} \label{def:bridge}  $K \subset S^3$ is in {\em
bridge position} (with respect to the height function $h$) if it is
in general position with respect to $h$ and all the minima of $K$
occur below all the maxima. A {\em minimal bridge
position} is a bridge position which minimizes the
number of minima and maxima. A level sphere lying between the minima
and maxima is called a {\em middle sphere}.
\end{definition}

Two bridge positions of a knot $K$ are regarded as equivalent if
they are isotopic via an isotopy throughout which the knot
remains in bridge position.  A knot may have more than one minimal
bridge position (see for example, Figure 13).  It is a tedious but
worthy exercise to show that if there is an isotopy between two knots
that are in bridge position and have the same middle sphere, and
during the isotopy the knot remains always transverse to the middle
sphere, then the bridge positions are equivalent.  In particular, a
sphere that divides a knot $K$ into two untangles, well
defines a bridge position of the knot in which the sphere is a
middle sphere.

\noindent{\bf Notation}\qua $\eta(X)$  means a regular neighborhood  of
$X$ in $Y$ for polyhedral spaces $X$ and $Y$.  For $K$ a knot in
the 3--sphere, put $E(K)=S^3-{\rm int}\eta(K)$.

Placing a knot in thin position can sometimes reveal essential
meridional planar surfaces in the complement, ie, planar
surfaces which are incompressible, have boundary a
collection of meridians of $K$, and which are not just annuli
parallel into the knot.  For example, there is this result of Thompson:

\begin{theorem}{\rm\cite{T}}\qua  If $E(K)$ does not have any essential
meridional  planar surface, then any thin position of $K$ is a minimal
bridge position of $K$ and vice versa.
\end{theorem}

In this paper we will be interested only in knots $K$ which have
tunnel number $1$.  That is, knots $K$ for which there is an arc
$\gam
\subset S^3$ such that $K \cap \gam = \bdd \gam$ and $S^3 - {\rm
int}\eta(K\cup\gamma)$ is a genus $2$ handlebody.  The arc $\gamma$ is
called an {\it unknotting tunnel} or {\it tunnel} for $K$.

Gordon and Reid have shown:

\begin{theorem}{\rm\cite{GR}}\qua Let $K \subset S^3 $ have tunnel number $1$.
Then  $E(K)$ does not have an essential meridional planar surface.
\end{theorem}

Combining these results gives the obvious corollary:

\begin{corollary}\label{cor:thin-bridge} Suppose $K$ has tunnel
number $1$.  Then any thin position of $K$ is a minimal
bridge position of $K$ and vice versa.
\end{corollary}

It is easy to see that all $2$--bridge knots have tunnel number one.
Moreover, Kobayashi has shown that the only tunnels possible are the
obvious ones:

\begin{theorem}{\rm\cite{K}}\qua  Any unknotting tunnel for a 2--bridge knot
is one of 6 known types.
\end{theorem}

These six types have the property that they can each be made level
with respect to the height function that gives the knot its $2$--bridge
structure.  That is, each can be put into a level sphere.  Morimoto
had earlier pointed out that this classification would follow rather
quickly if we knew that each tunnel could be made level:

\begin{theorem}{\rm\cite{M1}}\qua  If an unknotting tunnel for a 2--bridge
knot  can be put into a level sphere, then it is one of the 6 known
types.
\end{theorem}

So it is natural to turn Morimoto's theorem into a question: could it
be true that any unknotting tunnel for any tunnel number $1$ knot $K$
can be put into a level sphere, once $K$ is put in minimal bridge
position?  Here we allow the tunnel to be ``slid'' as well as
isotoped.  That is, we allow the ends of the tunnel to be moved
around on $K$,
eg, possibly past each other on $K$, or, indeed, one end of the tunnel may be
moved up onto the other end of the tunnel, changing what was an edge into
a  $1$--complex we will call an ``eyeglass'': an edge attached to a
circle at one of its ends.  Of course the reverse move is also
allowed.

In this paper we show that the answer to Morimoto's question is yes:

\begin{theorem} \label{theorem:main}  If $K \subset S^3$ is a tunnel
number one knot in minimal bridge position and $\gam$ is a tunnel
for $K$, then $\gam$ may be slid and isotoped to lie entirely in a
level sphere for $K$.
\end{theorem}

Moreover, in Section 6, we characterize the position of $\gam$
on a minimal bridge sphere.

This conclusion allows us to say a bit about the bridge structure of
tunnel number $1$ knots.

For example, \cite{GOT} gives a detailed description of tunnel
number $1$ links with the property that one component is the unknot.
The authors of \cite{GOT} ask [Problem 6.4] for an explicit 
description of tunnel
number one links both of whose components are knotted. In \cite{Be} Berge
constructed some examples by edge sliding, beginning with the case of
an unknotted component.

The theorem here assures us that such links can all
be obtained either in this way, or by beginning with a level
unknotting tunnel for a tunnel number one knot and sliding the ends
together by a (possibly complicated) path on the boundary of the knot
neighborhood. As a kind of generalization of these ideas we get the
following way of characterizing a tunnel number one knot as the
result of a series of bridge-increasing slides:

\begin{lemma} Suppose $\gam$ is a tunnel for a tunnel number $1$
knot $K$.  Let $m$ be a meridian of $\gamma$ on
$\partial\eta(K\cup\gamma)$. Suppose $K'$ is a simple closed curve on
$\partial\eta(K\cup\gamma)$ which intersects $m$ transversely once.
Then $K'$ is a tunnel number $1$ knot.
\end{lemma}

\begin{proof} Slide the ends of $\gamma$ around to ``peel off" $K'$.
\end{proof}

\begin{definition} We then say that $K'$ is obtained by a {\em tunnel
move} from $K$.
\end{definition}

Once we have established that such tunnels can always be made level,
we have:

\begin{proposition} Let $K$ be a tunnel number $1$ knot.  Then there
exists a sequence of knots beginning with the unknot $=K_0, K_1,
K_2,\ldots, K_n=K$ of bridge number $b_0 = 1, b_1, \ldots , b_n$ such
that
\begin{enumerate}
\renewcommand{\labelenumi}{\rm(\alph{enumi})}
\item
$K_i,\,i=1,\ldots ,n,$ is a tunnel number one knot,
\item
$K_i$ is obtained from $K_{i-1}$ by a tunnel move,
\item $b_i \geq 2 b_{i-1} - 1$.
\end{enumerate}
\end{proposition}

\begin{proof} By induction on bridge number, using the fact that
tunnels can always be made level by slides.
\end{proof}

The core of these results extends easily to tunnel number one links
and there is a brief discussion of this in the last section.

\section{Thin position and unknotting tunnels}

Let $K$ be a tunnel number $1$ knot and $\gamma$ be an unknotting
tunnel for $K$.  It will be convenient to think about the tunnel
$\gamma$ in two different ways.  Initially it is an arc with each of
its ends on $K$.  But the thinning process we will describe below may
force the ends of the arc together on $K$ until they are incident to
the same point on $K$ (so that $\gamma$ becomes a loop) and then
resolve the resulting $4$--valent vertex into two $3$--valent vertices
by pinching the ends of the tunnel together into a single arc.  (This
process could also be described as sliding one end of the tunnel onto
a neighborhood of the other end.)  The result is to change $\gamma$
into the union of an arc and a circle; as noted above, we will call
this graph an {\it eyeglass} and continue to denote it by $\gamma$.
No matter whether $\gamma$ denotes an edge or an eyeglass, the
neighborhood $\eta(K\cup\gamma)$ is a genus 2 handlebody $V_1$, and
its closed exterior $E(K\cup\gamma)$ is also a genus 2 handlebody
$V_2$.  Let $F$ denote the genus 2 surface that is their common
boundary, so $F$ is a Heegaard splitting surface for $E(K)$.

In this section we will extend Gabai's notion of thin position (see
Definition \ref{def:knotthin}) to the graph $K \cup \gamma.$  As
before, choose a height function
$h\co S^3-\{x,y\}=S^2\times \mathbb R \rightarrow\mathbb R$ and let
$S(t)=h^{-1}(t)$.

\begin{definition}
$K\cup\gamma$ is in {\it Morse position} with respect to $h$ if
\begin{enumerate}
\renewcommand{\labelenumi}{(\alph{enumi})}
\item the critical points of $h|K$ or $h|\gamma$ are nondegenerate
and those in $\gam$ lie in the interior of $\gamma$,
\item the critical points of $h|K$, $h|\gamma$ and the two vertices
in $K \cup \gamma$ all occur at different heights.
\end{enumerate}
\end{definition}

The heights at which there is a critical point  of either
$h|K$ or $h|\gamma$ or there is a vertex of $K \cup \gam$ are called
the {\it critical heights} for $K\cup\gamma$.  Near its
critical height, each vertex
$v$ can be classified into one of four types:
\medskip
\begin{enumerate}
\item
all ends of incident edges lie below $v$,
\item
  all ends of incident edges lie above $v$,
\item
exactly two ends of incident edges lie above $v$,
\item
exactly two ends of incident edges lie below $v$.
\end{enumerate}
\medskip
We will further simplify the local picture by isotoping  a small
neighborhood of a vertex of type 1 (respectively type 2),  transforming it
into a vertex of type 4 (respectively type 3) and a nearby maximum
(respectively minimum). The end result is called {\it normal form} for
$K\cup\gamma$.

\begin{definition}  A vertex of type 3 is called a {\em $Y$--vertex}
and a vertex of type 4 is called a {\em $\lambda$--vertex}. When we
refer to the {\em maxima} of $K \cup \gamma$ we will include all local
maxima of
$K$, all local maxima of $\inter (\gam) $, and all
$\lambda$--vertices.  Similarly, by the {\em minima} of $K \cup
\gamma$ we mean all local minima of
$K$, all local minima of $\inter (\gam)$, and all $Y$--vertices. A
maximum (respectively minimum) that is not a $\lambda$--vertex (respectively
$Y$--vertex) will be called a {\em regular} maximum (respectively minimum).
The union of the maxima and minima (hence including the vertices) are
called the {\em critical points} of $K \cup \gamma$ and their heights
the {\em critical values} or {\em critical heights}.
\end{definition}

\begin{definition} Extending and modifying the idea of knot width
$W(K)$ (see Definition \ref{def:knotthin}), let
$t_0 <
\ldots < t_n$ be the successive critical heights of $K \cup \gamma$
and suppose
$t_j$ and $t_k$ are the two levels at which the vertices occur. Let
$s_i, 1
\leq i \leq n$ be generic levels chosen so that $t_{i-1} < s_i <
t_i$.  Define $$W(K
\cup \gamma) = 2(\Sigma_{i \neq j,k}|S(s_i)\cap (K
\cup \gamma)|) +  |S(s_j)\cap (K \cup \gamma)| + |S(s_k)\cap (K \cup
\gamma)|.$$ \end{definition}

\begin{definition}
For $K$ in minimal bridge (hence thin) position, a {\it thin position}
of a pair $(K,\gamma)$ (rel $K$) is a position  which minimizes the
width of $(K,\gamma)$ (without changing the bridge presentation of
$K$).
\end{definition}

\begin{remark}  Of course, we can combine the two types of thinness:
Define the width $W(K,\gamma)$ of $(K,\gamma)$ to be the pair of
integers $(W(K),W(K \cup \gamma))$ in lexicographical order and
define a thin presentation of $(K, \gamma)$ to be one that minimizes
$W(K,\gamma)$.  This guarantees that the knot $K$ will
be thin. Note that there is no way of
predicting, {\em prima facie}, whether in a thin presentation $\gam$
will be an edge or an eyeglass.
\end{remark}

\begin{remark}  There is nothing magical about the description of
width above; thinning with respect to it serves the ad hoc purpose of
pushing maxima below minima, while remaining indifferent to pushing
maxima past maxima and minima past minima.  See the next section.
Other definitions that serve the same purpose are possible.  For
example, here's an alternative which might strike some as more
natural:  make
$V_1$ very thin and count components of intersection of level
spheres with
$V_1$ between successive critical points of $h$ on $\bdd V_1$.  Thus
two level spheres are used near a regular minimum or maximum of $K
\cup
\gamma$ but only one near a minimum or maximum that occurs at a
vertex. This definition of width is not equivalent to ours, but has
the same general properties that we need, as described in the next
section.\end{remark}

  Let $D$ be a meridian disk of $V_2 =
S^3-\inter(\eta(K\cup\gamma))$.

\begin{definition} \label{def:normal} The disk $(D,\partial D)\subset
(V_2, \bdd V_2)$ is in {\it
normal form} if
\begin{enumerate}
\renewcommand{\labelenumi}{(\alph{enumi})}
\item $\bdd D$ intersects each meridian of each edge of the graph $K
\cup \gamma$ minimally, up to isotopy in $\bdd V_2 = \bdd\eta(K
\cup \gam)$,
\item  each critical point of $h$ on $D$ is nondegenerate,
\item  no critical point of $h$ on int$(D)$ occurs at a critical
height of
$K\cup\gamma$,
\item  no two critical points of $h$ on int$(D)$ occur at the same
height,
\item  the minima (respectively maxima) of $h|\partial D$ at the minima
(respectively maxima) of $K \cup \gamma$ are also local extrema of $h$ on
$D$, ie,  `half-center' singularities,
\item the maxima of $h|\partial D$ at $Y$--vertices and the minima of
$h|\partial D$  at $\lambda$--vertices are, on the contrary,
`half-saddle'  singularities of $h$ on $D$.
\end{enumerate}
\end{definition}

Standard Morse theory ensures that, for $K \cup \gamma$ in normal
form, any properly imbedded essential disk  $(D,\partial D) \subset
(V_2, \bdd V_2)$ can be put
in normal form.

\begin{definition}
$K\cup\gamma$ is in {\it bridge position} if  there is no
minimum of $K\cup\gamma$ at a height above a maximum.
\end{definition}

\begin{definition} Let $S$ be a generic level 2--sphere and  $B_u$ and
$B_l$ denote the balls which are the closures of the region above
$S$ and below $S$ respectively. An {\it upper disk} (respectively {\it
lower disk}) for $S$ is a disk $F \subset V_2$  transverse to
$S$ such that $\partial F=\alpha\cup\beta$, where $\alpha = \bdd F
\cap S$ is  an arc properly imbedded in $S - \eta(K \cup \gamma)$,
$\beta$ is an arc imbedded  on
$\partial\eta(K
\cup
\gamma)$, $\partial\alpha=\partial\beta$ and a small product
neighborhood of $\bdd F$ in $F$ lies in $B_u$ (respectively $B_l$)  ie, it
lies {\it above} (respectively {\it below}) $S$.
\end{definition}

Note that $\inter(F)$ may intersect $S$ in  simple closed curves.

\section{Moves that thin and moves that don't}\label{section:thinning}

Width is defined as it is in order to ensure that certain types of
moves from one normal form of $K \cup \gamma$ to another will
decrease width, some will increase width, and others will have no
effect.  We list the possible ways in which adjacent minima and maxima
can be moved: Suppose, as in the definition of $W(K \cup \gamma)$,
$s_i$ is a regular level, lying between adjacent critical levels or
heights of vertices $t_{i-1}$ and $t_i$, with $t_{i-1} < s_i < t_i$.
For each $1 \leq i \leq n,$ let $p_i = |S(s_i) \cap (K \cup \gamma)|.$

We will assume there is an isotopy of $K \cup \gamma$ whose only
effect on the positioning of critical points is to exchange the
critical points at two adjacent levels, say
$t_{i-1}$ and $t_i$. (We say that the critical point at
$t_{i-1}$ is {\em moved} past the critical point at $t_i$.)  For this
to make easy sense, we must assume that no subarc connecting the two
critical points lies completely between the two levels.

Simple computations (see Figure 1) show the following effects:

\begin{claim} If the critical points at $t_{i-1}$ and $t_i$ are both
maxima or both minima, moving one past the other has no effect on the
width.
\end{claim}

\begin{proof}
This is obvious if they are both regular maxima, both regular minima,
both $Y$--vertices or both $\lambda$--vertices.  The remaining
representative case is this: at $t_i$ there is a $Y$--vertex and at
$t_{i-1}$ there is a regular minimum.

Only $p_i$ is changed (to $p'_i = p_i - 1$).  But because the
height of the $Y$--vertex changes, note that, before the change, the
term $2p_{i-1} + p_i$ appears in the definition of $W(K \cup \gamma)$
and afterwards the term $p_{i-1} + 2p'_i$ appears.  But since
$p_{i-1} = p_i - 2 = p'_i-1$ the change has no effect.
\end{proof}

\begin{figure}[ht!]
\cl{\small
\SetLabels 
(0.2*-0.1){$2p_{i-1}+p_i=5$}\\
(0.8*-0.1){$p_{i-1}+2p_i=5$}\\
(0.5*0.5)$s_i$\\
(0.5*0.2)$s_{i-1}$\\
\endSetLabels 
\AffixLabels{\includegraphics[width=.9\textwidth]{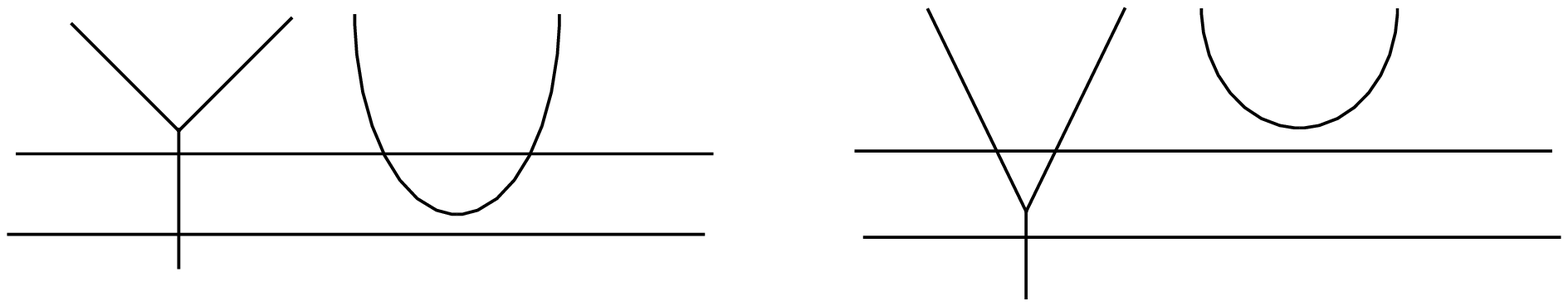}}}
\caption{}
\end{figure}

\begin{claim} If the critical point at $t_{i-1}$ is a minimum
and that at $t_i$ is a maximum then $W(K \cup \gam)$ is reduced.
\end{claim}

\begin{proof}
Only $p_i$ is affected. If both critical points are regular critical
points, $p_i$ is reduced by $4$ so $W(K\cup\gamma)$ is reduced by
$8$.

If both are vertices, then $p_i$ is reduced by $2$, so
$W(K\cup\gamma)$ is reduced by $2$.

If the critical point at $t_{i-1}$ is a
$Y$--vertex and the one at $t_i$ is a regular maximum or if  the
critical point at $t_{i-1}$ is a  regular minimum and the one at
$t_{i}$ is a $\lambda$--vertex, then $p_i$ is reduced to $p'_i = p_i -
3$. The effect on the calculation of $W(K\cup\gamma)$ is, in the
first case, to replace the term $p_{i-1} + 2p_i$ with $2p_{i-1} +
p'_i$.  Since $p_{i-1} = p_i - 1$ the net effect is to reduce
$W(K\cup\gamma)$ by $4$. (See Figure 2.) In the second case the term
$2p_{i-1} + p_i$ is replaced by $p_{i-1} + 2p'_i$.  Since here
$p_{i-1} = p_i - 2$ the net effect is again to reduce
$W(K\cup\gamma)$ by $4$.\end{proof}

\begin{figure}[ht!]
\cl{\small
\SetLabels 
(0.2*-0.1){$p_{i-1}+2p_i=11$}\\
(0.8*-0.1){$2p_{i-1}+p_i=7$}\\
(0.5*0.5)$s_i$\\
(0.5*0.2)$s_{i-1}$\\
\endSetLabels 
\AffixLabels{\includegraphics[width=.9\textwidth]{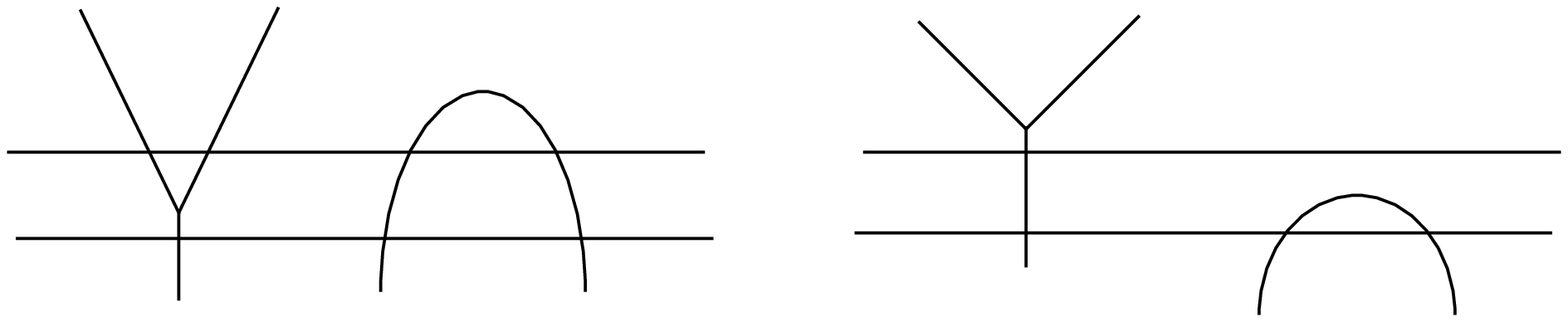}}}
\caption{}
\end{figure}

Now we consider the case in which the critical points at levels
$t_{i-1}$ and $t_i$ are {\em adjacent}.  That is, suppose there is a
subarc of $K \cup \gamma$ that descends from one to the other.  If
one is a maximum and the other a minimum, it may be possible to

\begin{itemize}

\item cancel the pair (if neither is a vertex),

\item cancel a regular maximum (respectively minimum) by changing a
$Y$--vertex (respectively $\lambda$--vertex) to a $\lambda$--vertex (respectively
$Y$--vertex), or

\item replace a $Y$--vertex at $t_{i-1}$ and $\lambda$--vertex at $t_i$
with a $\lambda$--vertex at $t_{i-1}$ and $Y$--vertex at $t_i$.

\end{itemize}

It's easy to check that in all these cases, $W(K \cup \gamma)$ is
reduced (or, in the last case, it is possible that $\gam$ can be
directly levelled).  The upshot of all these comments is the
following overarching principle:

\begin{proposition} \label{prop:overarch}

Suppose $K$ is in minimal bridge position and $S$ is a generic level
sphere lying just above a maximum of
$K \cup \gamma$.  Suppose a subarc $\alpha \subset (K \cup \gamma)$
doesn't contain the maximum, has ends incident to $S$ from below, and
$\alpha$ can be moved (by slide and/or isotopy)
to lie on $S$, with its ends incident to the ends of $(K \cup
\gamma) - \alpha$, as before.  Then this move will reduce $W(K \cup \gam)$.  A
symmetric statement holds for $S$ just below a minimum. (See Figure
3.)

\end{proposition}

\begin{proof}
Conceptually, in terms of the height function on $\alpha$, the arc
$\alpha$ can be moved, leaving its complement $(K \cup \gam) - \alpha$
fixed, in a process that only moves maxima down, minima up, with
possible cancellations occurring.  (There is no claim that the actual
motion of $\alpha$ can be made to realize this conceptual process.)
These moves can only thin and can never thicken.  The arc $\alpha$ has
at least one minimum, and it is eventually moved above the maximum
that lies just below $S$.  Hence, at least by this move, the
presentation is thinned.
\end{proof}

\begin{figure}[ht!]
\cl{\small
\SetLabels 
(0.15*.5)$\alpha$\\
(0.41*0.80)$S$\\
(0.97*0.80)$S$\\
(0.5*0.63)$\scriptstyle\rm{Thin}$\\
\endSetLabels 
\AffixLabels{\includegraphics[width=.8\textwidth]{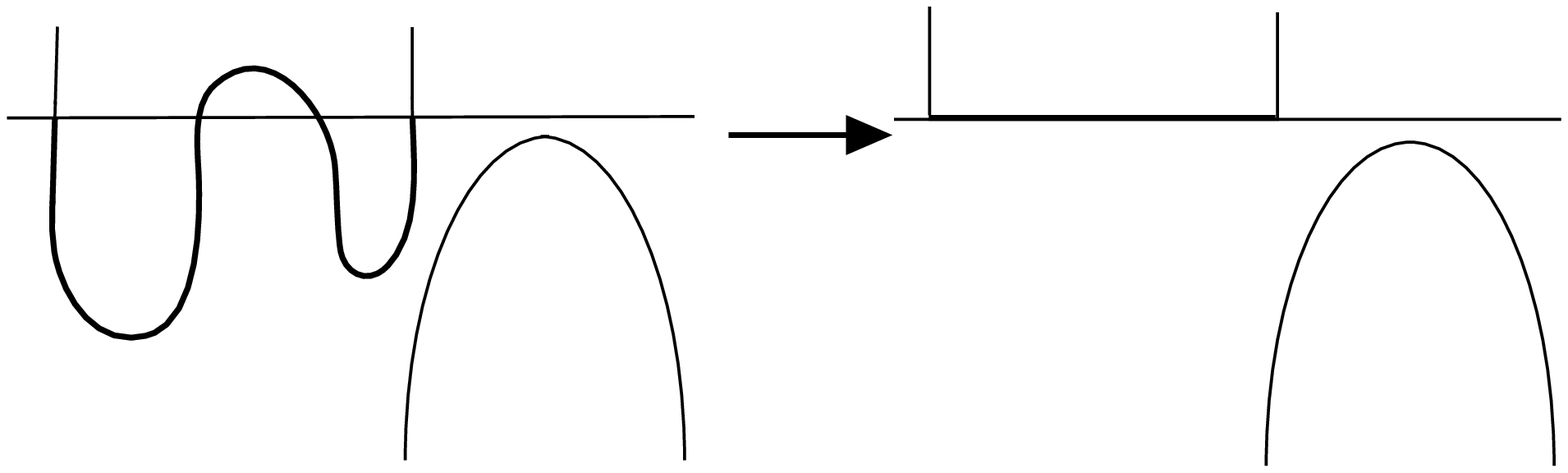}}}
\caption{}
\end{figure}

Note that the move described in Proposition \ref{prop:overarch} could
change the bridge presentation of $K$.  This can happen only if $\aaa$
has a single minimum (because otherwise $K$ itself would have been
thinned by the process) and, during the isotopy of $\aaa$,
there are times when $\aaa$ has more critical points than a single
minimum.

\section{A thin $K \cup \gamma$ is in bridge position}
\label{section:thin=bridge}

Just as Corollary \ref{cor:thin-bridge} shows that if $K$ is in thin
position it is  in bridge position, we will show in this section
that if $K\cup\gamma$ is put in thin position rel $K$ it is  also in
bridge position.

Recall that the neighborhood
$V_1 = \eta(K\cup\gamma)$ is a genus 2 handlebody,
as is its closed exterior $V_2 = E(K\cup\gamma)$.  We begin by
recalling some crucial facts about surfaces in genus two handlebodies.

\begin{definition} Let $M$ be a compact 3--manifold and $P$ a compact
2--manifold  properly imbedded in $M$. We say that $P$ is {\it
$\partial$--parallel} if $P$ is isotopic rel $\partial P$ to a compact
manifold in $\partial M$,  and that $P$ is {\it essential} if $P$ is
incompressible and has  a component which is not $\partial$--parallel.
\end{definition}

\begin{proposition}{\rm\cite{M2}}\qua  \label{prop:morimoto} Let $M$ be an
orientable closed 3--manifold with a genus 2  Heegaard splitting
$(V_1,V_2)$. If $M$ contains a 2--sphere
$S$ such that each component of $S\cap V_1$ is a non-separating  disk
in $V_1$ and $S\cap V_2$ is an essential planar surface  in $V_2$,
then $M$ has a lens space or $S^2\times S^1$ summand.
\end{proposition}

\begin{proposition}{\rm\cite{M2}}\qua  \label{prop:morimoto2} Let $V$ be a
genus $g>1$ handlebody, and $P$ be a connected incompressible  planar surface
with $\ell>1$ boundary components properly imbedded in $V$. If
$\partial P$ consists of mutually parallel separating loops in
$\partial V$, then $\ell=2$ and $P$ is a $\bdd$--parallel annulus.
\end{proposition}

(In fact, Morimoto also shows that when $g = 2$ we can drop the hypothesis
that the loops are separating.)

\bigskip

We proceed then to prove:

\begin{proposition}\label{prop:thin=bridge} Suppose that
$K$ is in minimal bridge position and $(K,\gamma)$ is in thin
position rel $K$. Then
$K\cup\gamma$ is in bridge position.
\end{proposition}

\begin{proof} If not, there is a regular value $t_0$ such that
$S(t_0)$ lies  between adjacent critical values $x$ and $y$ of $h$,
where
$x$ is a minimum of $K \cup \gamma$ lying above $t_0$ and $y$ is a
maximum  of
$K \cup \gamma$ lying below $t_0$. Let $P$ be the planar surface
$S(t_0)-$int$(\eta(K\cup\gamma))$.

\medskip

\noindent {\bf Case A}\qua $\gamma$ is a single edge (ie, $\gam$ is not an
eyeglass).

Compress $P$ as much as possible in $E(K\cup\gamma)$.  The
compressions may take place to either side of $P$, and may  have to
be done in several steps. Let $\tilde{P}$ be the resulting
collection of meridional planar surfaces. $\tilde{P}$ is
incompressible  in $E(K\cup\gam)$. Since $\gamma$ is an edge, each component
of $\partial\tilde{P}$ is non-separating in
$\partial\eta(K\cup\gamma)$.  Thus, by Proposition
\ref{prop:morimoto}, each component of $\tilde{P}$ is boundary
parallel. Let $P'$ be a component of $\tilde{P}$, so
$S'=P'\cup$(meridian disks of $K\cup\gam$) bounds a 3--ball, one of 5 cases as
in Figure 4. In each case, we can find a disk $E$  as in Figure 4.
That is,  $\partial E$ consists of two arcs $\alpha$  and $\beta,
\partial\alpha=\partial\beta$, with $\alpha$ a properly imbedded
essential arc in $P'$, and $\beta$ an arc imbedded on $\partial
E(K\cup\gamma)$. Moreover,
\begin{enumerate}
\item
$\beta$ is parallel to an arc of $K$, ($P'$ is an annulus),
\item
$\beta$ is parallel to an arc of $\gamma$, ($P'$ is an annulus),
\item
$\beta$ is parallel to an arc of $K\cup\gamma$  which contains a
vertex, ($P'$ is a three punctured sphere),
\item
$\beta$ is parallel to an arc of $K\cup\gamma$  which contains a
vertex, ($P'$ is a four punctured sphere),
\item
$\beta$ is parallel to whole $\gamma\cup$(two arcs) of $K$,  ($P'$ is
a four punctured sphere).
\end{enumerate}

\begin{figure}[ht!]
\cl{\small
\SetLabels 
(0.1*0.96)$E$\\
(0*0.84)$\alpha$\\
(0.22*0.80)$K$\\
(0.305*0.90)$P'$\\
(0.65*0.91)$\gamma$\\
(0.9*0.79)Case 1\\
(0.9*0.59)Case 2\\
(0.9*0.39)Case 3\\
(0.9*0.20)Case 4\\
(0.9*-0.02)Case 5\\
\endSetLabels 
\AffixLabels{\includegraphics[width=.35\textwidth]{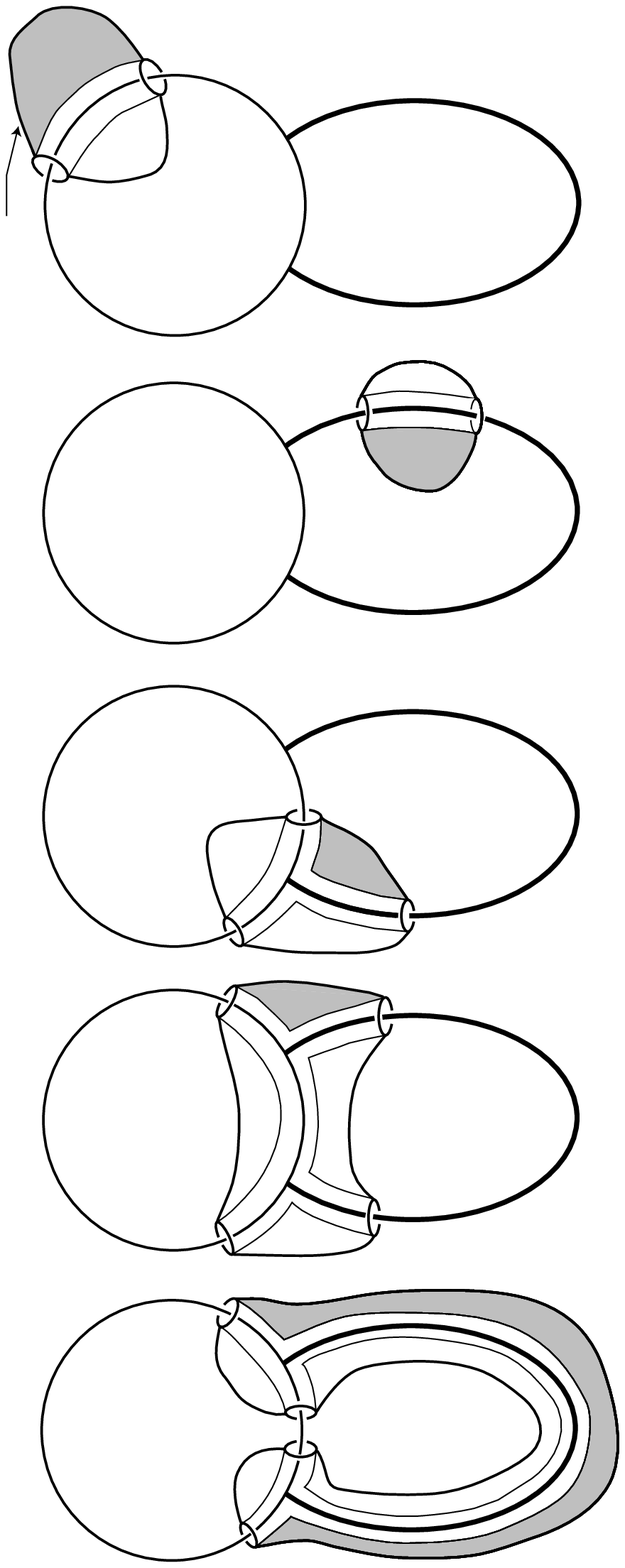}}}
\caption{}
\end{figure}

Now reverse the compressions  on $\tilde{P}$ to reassemble $P$. This
means that the components of
$\tilde{P}$ are reattached by tubes, which may run through each
other,  and which may intersect $E$. However, since the tubes can be
chosen  to be disjoint from a neighborhood of $\partial E$,
$\aaa$ remains an arc entirely on $S(t_0)$ and the ends of $\bbb$ lie
either both below or both above $S(t_0)$, say below.  Then by a slide
or isotopy $E$ can be used to move $\bbb$ to $\aaa$, and thus above
level $y$.  This contradicts thinness, via eg, Proposition
\ref{prop:overarch}. (In particular, in Case
5,  $\gamma$ becomes level and we are done.)

All but the first case can be viewed as slides and isotopies of
$\gam$ and so do not affect the bridge presentation of $K$.  So to
finish Case A it suffices to show that not all components of
$\tilde{P}$ are boundary parallel annuli.  The first observation is
that if any of these annuli were nested then the thinning move on
$E$, defined for an annulus that is not innermost, would actually
thin $K$, which we know to be impossible.  Hence the annuli are not
nested.  For the same reason we know that the segments of $K$ which
are cut off by the annuli must all lie above (or all below) the sphere
$S(t_0)$.  So if these boundary parallel annuli constitute all of
$\tilde{P}$ then the only part of $K \cup \gam$ that lies above
$S(t_0)$ are the arcs of $K$ which the annuli cut off.  In
particular, the minimum of $K \cup \gam$ that lies just above
$S(t_0)$ must also belong to $K$.  But this is a contradiction, for
the isotopy of $K$ which moves the segment containing the minimum
down to $S(t_0)$ would thin $K$ itself.

\medskip

\noindent {\bf Case B}\qua $\gamma$ is an eyeglass.

We denote by $\gamma_{c}$ the circle part of $\gamma$,  and by
$\gamma_{a}$ the arc part of $\gamma$ as in Figure 5. If the
link $K \cup \gam_c$ is disjoint from $S(t_0)$, the proof proceeds
much as in Case A: Compress
$P$ as much as possible in the  complement of
$K \cup \gam$ and let $P'$ denote a component of the result.  Since
$P'$ intersects $K \cup \gamma$ only in $\gamma_a$, it follows from
Proposition \ref{prop:morimoto2} that $P'$ is a
$\bdd$--parallel annulus. This implies, as in Case A, that a subarc of
$\gamma_a$ can be isotoped to $P$, contradicting thinness via
Proposition \ref{prop:overarch}. $K$ is left unchanged.
An argument as in Case A-5 similarly applies if $\gam_{a}$
is disjoint from $S(t_0)$.

It remains to consider the case in which $K$ or $\gam_{c}$
and also $\gam_a$ do intersect $S(t_0)$; a preliminary
construction may be needed.
By a slide, make $\gamma$ into a
simple arc $\gamma'$. (Of course, after the slide
$K\cup\gamma'$ may no longer be in thin position.)  The new arc
$\gamma'$ consists of three subarcs: $\gamma_c$ (with a small subarc
$\gamma_{\epsilon}$ removed near its vertex),  $\gamma_a$, and a new
arc $\gamma'_a$ parallel to
$\gamma_{a}$.  Let $G$
be a disk  which gives the parallelism of $\gamma_a$ and $\gamma'_a$
(see Figure 5). The boundary of $G$ consists of four arcs:
$\gamma_a$,
$\gam'_a$, $\gamma_{\epsilon}$, and a short subarc $K_{\epsilon}$
of
$K$.  We may take $\gamma_{\epsilon}$ and
$K_{\epsilon}$ to be disjoint from $S(t_0)$.

As in Case A, compress $P$ as much as possible in the  complement of
$K \cup \gam'$ and let $P'$ denote a component of the result, chosen
to have at least one boundary component a meridian of either $K$ or
$\gamma_c$.  From Proposition \ref{prop:morimoto} $P'$ is
$\bdd$--parallel in $V_2=E(K\cup\gam ')$.
After the compressions, it may be that $P'$
intersects $\gamma_{\epsilon}$, but it must remain disjoint from
$K_{\epsilon}$. We have 4 cases as illustrated in Figure 6.

In each case, we can find a disk $E$  as in Figure 6.
That is,  $\partial E$ consists of two arcs $\alpha$  and $\beta,
\partial\alpha=\partial\beta$, with $\alpha$ a properly imbedded
essential arc in $P'$, and $\beta$ an arc imbedded on $\partial
E(K\cup\gamma')$.  Note that
$E \cap G=\emptyset$.  We write as in Case A:
\begin{enumerate}
\item
$\beta$ is parallel to an arc of $K$, ($P'$ is an annulus),
\item
$\beta$ is parallel to an arc of $\gamma_a\cup\gamma_c$ or
$\gamma'_a\cup\gamma_c$, ($P'$ is an annulus),
\item
$\beta$ is parallel to an arc of $\gamma_c$, ($P'$ is an annulus),
\item
$\beta$ is parallel to an arc of $K\cup\gamma_a$,  ($P'$ is a four
punctured sphere).
\end{enumerate}

Now reverse the compressions  on $\tilde{P}$ to
reassemble $P$. This means that the components of
$\tilde{P}$ are reattached by tubes, which may run through each
other,  and which may intersect $E$ and $G$. However,  $E$ and
$G$ remain disjoint, and $\bdd E$ remains unchanged since the tubes
can be chosen  to be disjoint from a neighborhood of $\partial E$.
So $E$ persists as a disk which can be used to move,  by isotopies
and edge slides, some subarc of the {\em original} $K \cup \gamma$
to the level sphere $S(t_0)$,  reducing $W(K \cup \gam)$ via
Proposition \ref{prop:overarch}.

As in Case A, this thinning can be done without affecting the bridge
position of $K$ unless all components of $\tilde{P}$ are boundary
parallel annuli; but this would lead to the same contradiction as in
Case A.
\end{proof}

\begin{figure}[ht!]
\cl{\small
\SetLabels 
\E(0.06*0.3)$K$\\
\E(0.22*0.65)$\gamma_a$\\
\E(0.79*-0.1)$\gamma'_a$\\
\E(0.34*0.3)$\gamma_c$\\
\E(0.65*0.5)$K_\epsilon$\\
\E(0.92*0.5)$\gamma_\epsilon$\\
\endSetLabels 
\AffixLabels{\includegraphics[width=.8\textwidth]{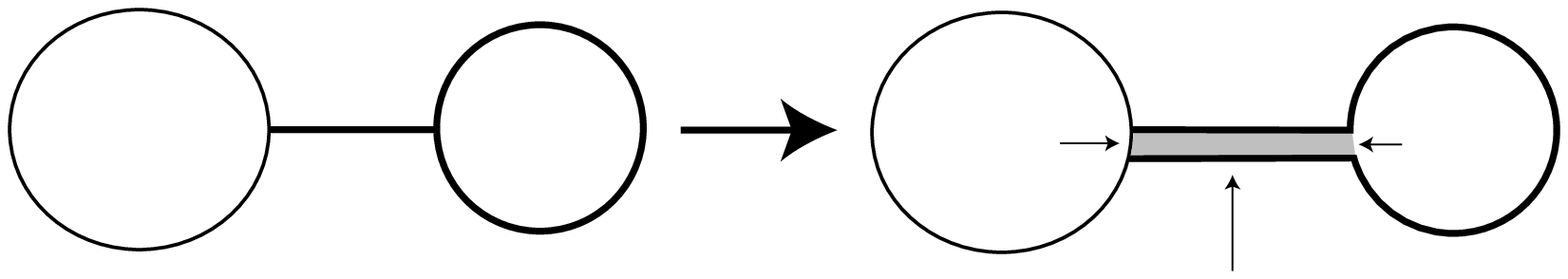}}}
\caption{}
\end{figure}

\begin{figure}[ht!]
\cl{\small
\SetLabels 
\E(0.29*0.55)Case 1\\
\E(0.8*0.55)Case 2\\
\E(0.29*-0.02)Case 3\\
\E(0.8*-0.02)Case 4\\
\E(0.16*0.65)$K$\\
\E(0.14*0.8)$P'$\\
\E(0.05*0.93)$E$\\
\E(0.29*0.75)$\scriptstyle G$\\
\B(0.29*0.8)$\gamma_a$\\
\T(0.29*0.7)$\gamma'_a$\\
\B(0.42*0.635)$\gamma_c$\\
\E(0.83*0.88)$E$\\
\E(0.43*0.16)$E$\\
\E(0.725*0.29)$\scriptstyle E$\\
\endSetLabels 
\AffixLabels{\includegraphics[width=.8\textwidth]{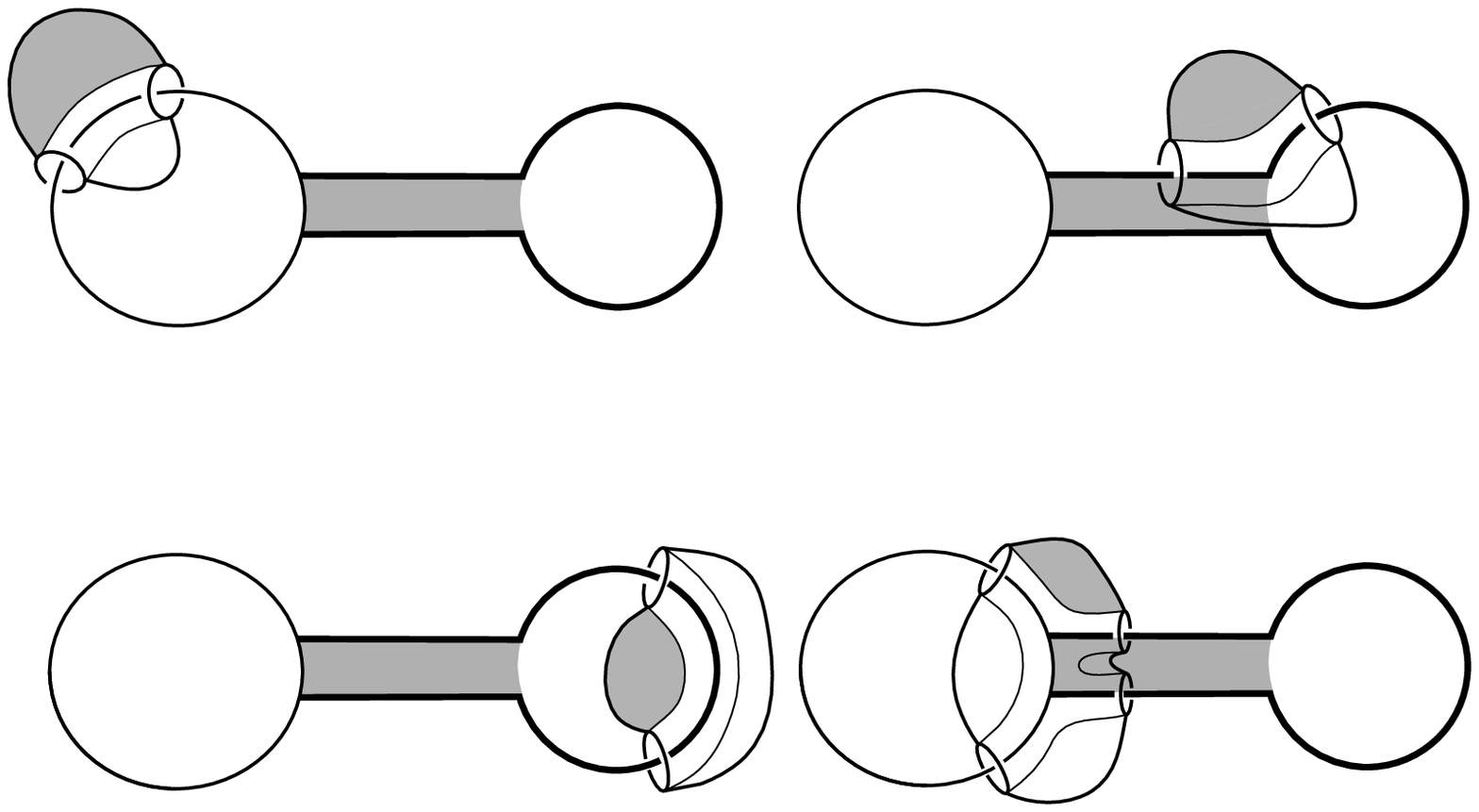}}}
\caption{}
\end{figure}

\section{Proof of Theorem \ref{theorem:main}}

Once $K \cup \gam$ is in bridge position, one can push any $Y$--vertex
below any distant minimum and any $\lam$--vertex above any distant
maximum without increasing the width or altering the bridge
structure of $K$.  So, typically, it is possible to choose a thin
position in which neither the highest local minimum nor the lowest
local maximum occurs at a vertex.  The single possible exception (for
which the classification of tunnels is anyway already known) is the
case in which $K$ is a
$2$--bridge knot and
$\gam$ is an arc with exactly one critical point in its interior (see
Figure 7).

\begin{figure}[ht!]
\cl{\small
\SetLabels 
\E(0.94*0.3)$K$\\
\B(0.55*0.6)$\gamma$\\
\endSetLabels 
\AffixLabels{\includegraphics[width=.4\textwidth]{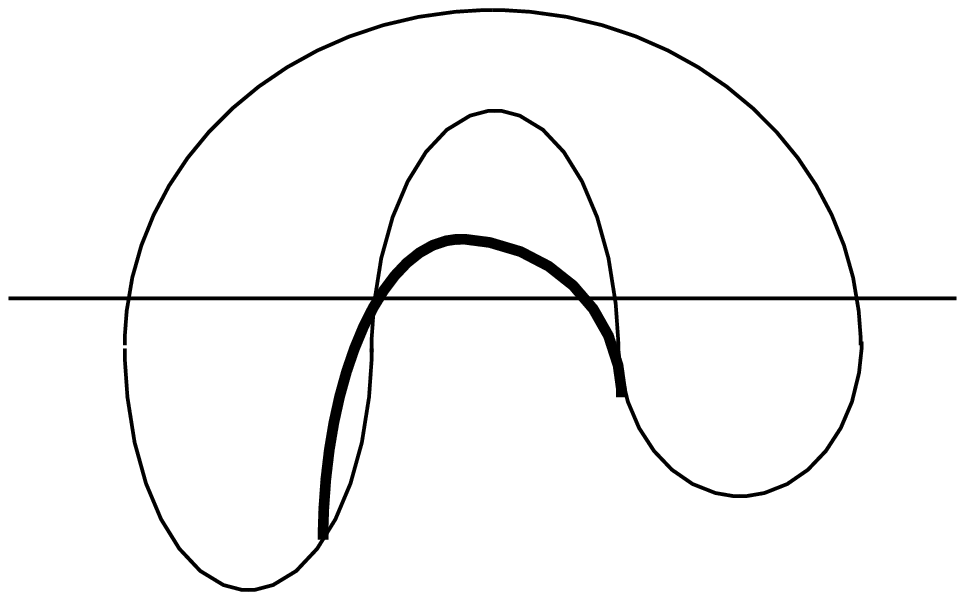}}}
\caption{}
\end{figure}

\begin{lemma} \label{lemma:level} Suppose $K$ is a tunnel number $1$
knot in minimal bridge, hence thin, position.  Then  $K \cup \gamma$
can be put in thin position rel $K$ so that, for $x$ the level of the
highest local minimum and  $y$ the level of the lowest local maximum,
there is a $z\in(x,y)$ so that $S(z)$ simultaneously has an upper
disk and a lower disk, and these disks are disjoint.
\end{lemma}

\begin{proof} This proof is essentially that of Lemma 4.4 in \cite{G}.
As above, let $D$ be a meridian disk for $V_2$. First observe that
$\bdd D$ runs along every edge of the graph
$K \cup \gam$, for otherwise it would follow that $K$ is
trivial.   As in Section 4 in
\cite{G}, by an isotopy,  we may assume that the local picture of $D$
near regular maxima and regular minima is as
illustrated  in Figure 8 (see also Definition
\ref{def:normal} (e)).

\begin{figure}[ht!]
\centering
\includegraphics[width=.35\textwidth]{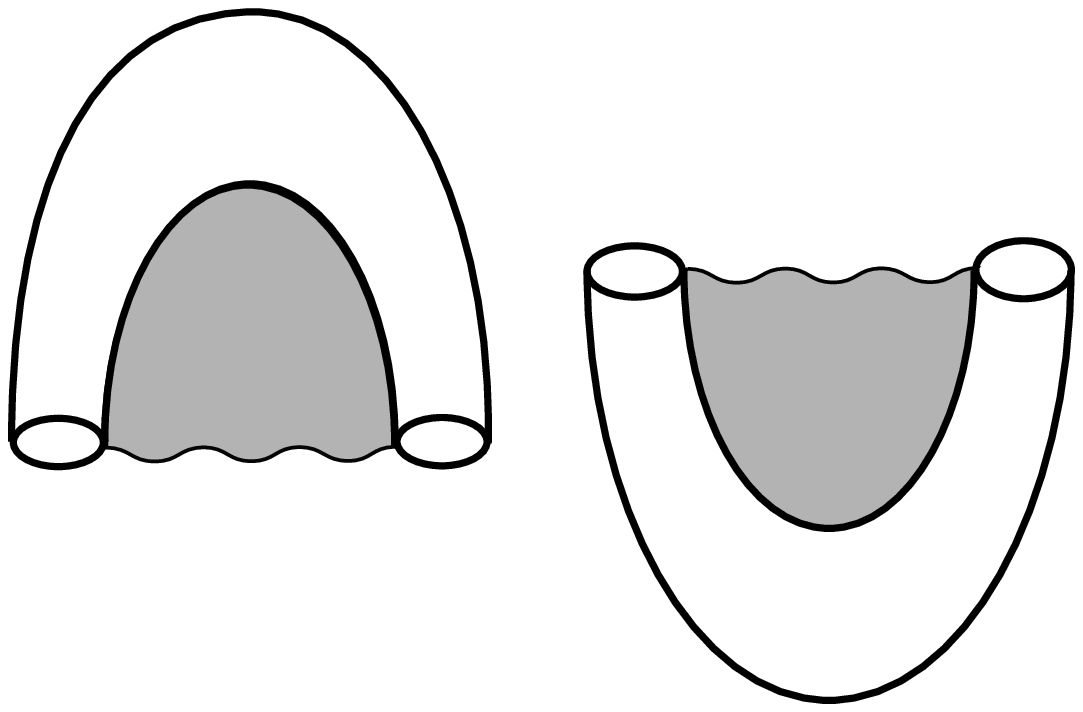}
\caption{}
\end{figure}

With the exception noted above (when $K$ is a
$2$--bridge knot and $\gam$ is an arc with exactly one critical point
in its interior), we may assume that the minimum at
$x$ and the maximum at $y$ are regular critical points, not vertices.
In particular, from Figure 8 we see that there is a small value
$\varepsilon$  such that $S(t)$ cuts off an upper outermost disk from
$D$ for  $t\in[y-\varepsilon,y]$ and cuts off a lower outermost disk
from $D$ for $t\in[x,x+\varepsilon]$.  At any regular level $t$
between $x$ and $y$ we have $(K \cup \gamma) \cap S(t) \neq
\emptyset$ so there is either an upper or a lower disk.  Since at the
lowest level there are lower disks and the upper level there are
upper disks and at every regular value there is one or the other, it
would seem to follow immediately that at some level between $x$ and
$y$ there are both upper and lower disks, as required.

But there remains the possibility that passing upward through a
critical level $t_c$ for $D$ (necessarily corresponding to a saddle
intersection of $S(t_c)$ with the interior of $D$), two lower disks
disappear while two upper disks are created. In this case, simply
thicken $D$ very slightly, creating a collar whose boundary consists
of two disjoint copies of $D$, each of them transverse to
$S(t_c)$. $S(t_c)$ will cut off from one an upper disk and
from the other a lower disk.

This completes the proof for all but the exceptional case: Suppose
there is a single critical point in the interior of $\gamma$, say a
maximum, and the ends of $\gamma$ are incident to $K$ just
above distinct minima of $K$, of which there are exactly two.
(See Figure 7.)  Then the highest minimum on $K \cup \gamma$ occurs
at a $Y$--vertex, from which one end of $\gamma$ ascends. The argument
above still applies, unless the local picture of $D$ near the
$Y$--vertex looks as in Figure 9, that is, unless the height
function on $D$ has one or more half-saddles at
$x$. This situation is characterized by the condition that
$\bdd D$ intersects a meridian of the edge descending from the
$Y$--vertex in more points than the sum of the number of points it
intersects the meridians ascending from the $Y$--vertex.

There are two possible
complications in this case:  One is the possibility that there could
be upper disks persisting just below level $x$.  We leave as an
exercise for the reader, most easily done after seeing the proof of
Theorem \ref{theorem:level1}, that such an upper disk could be used
to thin the presentation further, a contradiction.  A more subtle
problem is that an upper disk $D_u$ for $t$ slightly above $x$ could
be destroyed by the half-saddle singularities at $x$ and one of
the pieces become a lower disk $D_l$ for $t$ slightly below $x$.

\begin{figure}[ht!]
\centering
\includegraphics[width=.25\textwidth]{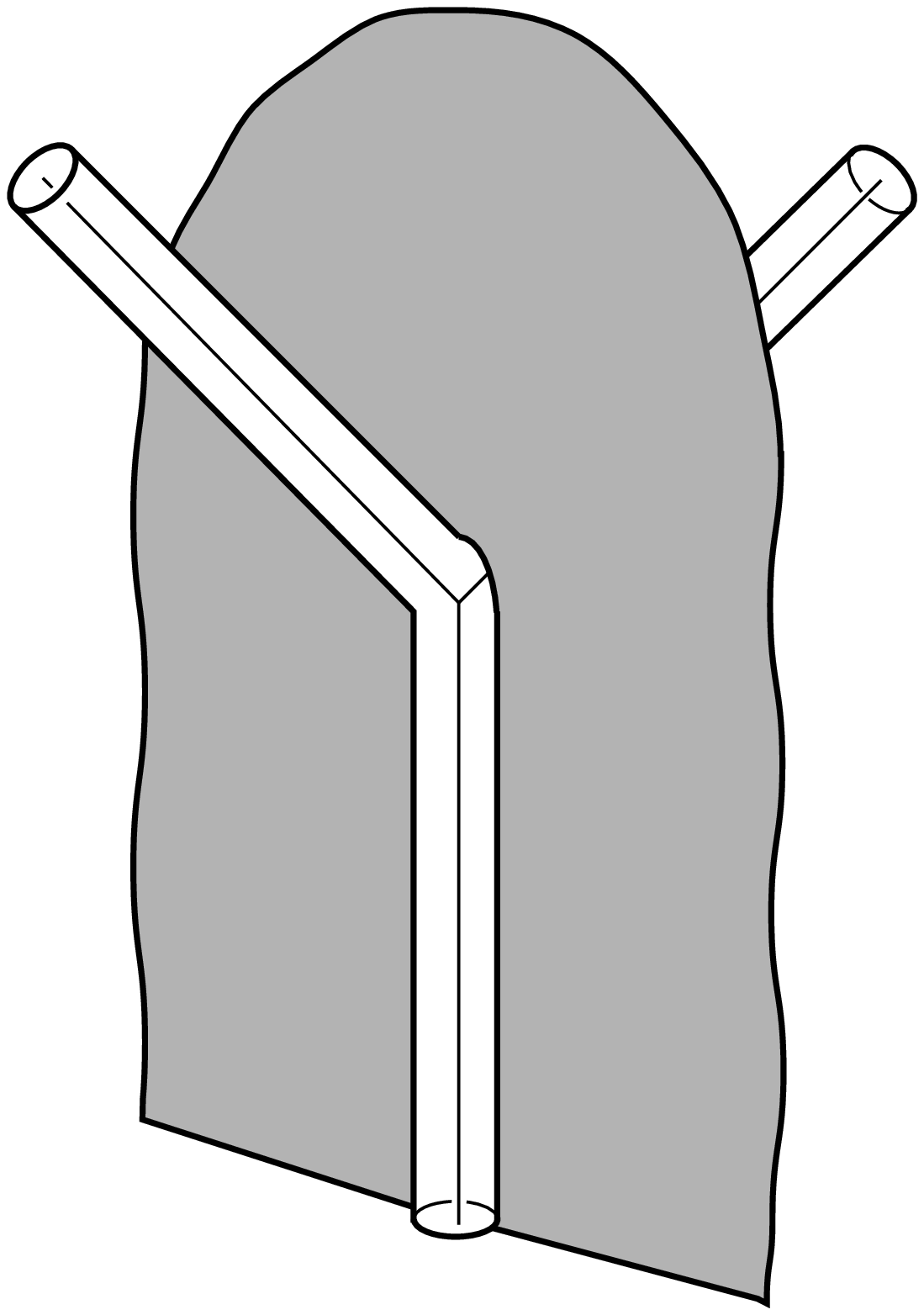}
\caption{}
\end{figure}

Note that the lower disk intersects $K \cup \gamma$
in a subarc of $K$, namely the arc around the minimum of $K$
to which the $Y$--vertex is adjacent.  This lower disk can be used to
slide the end of $\gamma$ at the $Y$--vertex across the minimum at $K$
and back to height $x$ so that afterwards no new critical points are
introduced on $\gamma$. (See Figure 10.) In particular, the placement
of $K \cup \gamma$ is still thin. But now the intersection of $D$ with
a neighborhood of the (new) $Y$--vertex is of a form that guarantees
a lower disk for $t\in(x,x+\varepsilon]$, since it is now the
meridian of an ascending arc from the $Y$--vertex that intersects
$\bdd D$ most often.  The argument then proceeds as above.
\end{proof}
\begin{figure}[ht!]
\cl{\small
\SetLabels 
\B(0.5*0.65)slide\\
\B(0.25*0.65)$\gamma$\\
\B(0.75*0.65)$\gamma$\\
\endSetLabels 
\AffixLabels{\includegraphics[width=.8\textwidth]{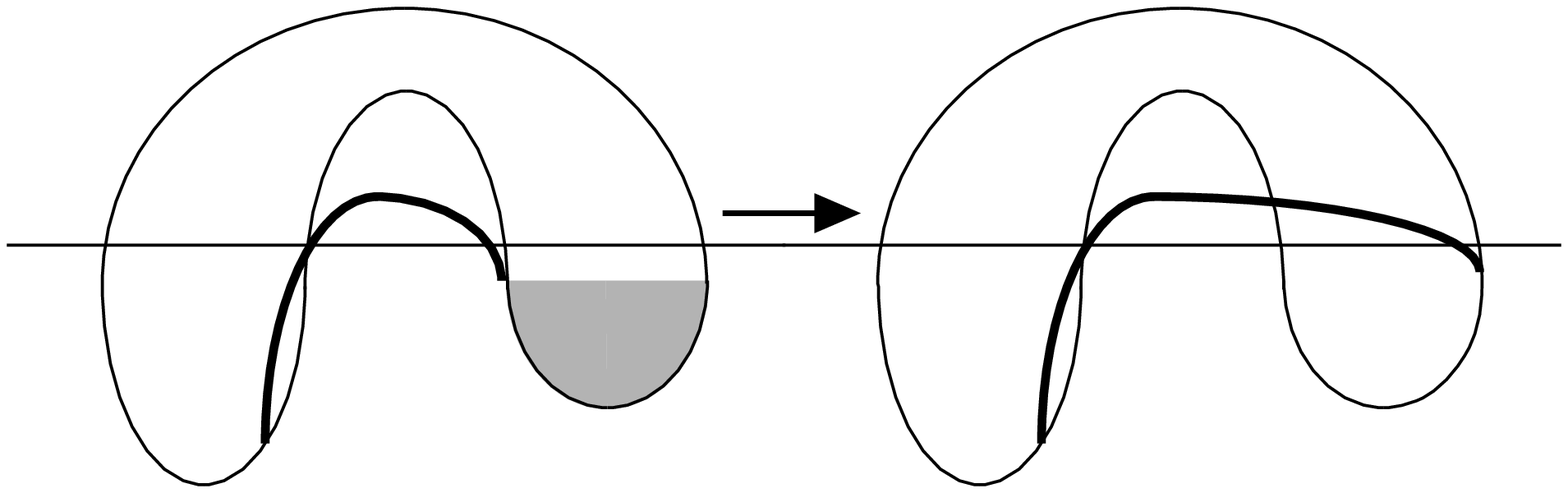}}}
\caption{}
\end{figure}

The previous proof shows that, even once $K \cup \gam$ is in thin
position, it may be useful to choose a particular thin position.
In fact, we will need the following highly technical lemma.

Suppose,
when $K \cup \gam$ is put in thin position, some edge
$e \neq \gamma$ of the graph $K \cup \gamma$ monotonically
descends from one $\lambda$--vertex to another. This is equivalent to
the condition that, for $S(z)$ a middle sphere, some edge of
$K \cup \gamma$ lies entirely above $S(z)$ or, equivalently, some
component $X$ of $(K\cup\gamma)-S(z)$ lying above $S(z)$ is a tree
with four ends and all other components of $(K\cup\gamma)-S(z)$ are
arcs.   Let $P$ be the $4$--punctured sphere $\eta(X) \cap \bdd \eta
(K \cup
\gam)$ and $C \subset P$ be a properly imbedded collection of
essential arcs, with the property that the $4$ components of $\bdd
P$ can be grouped into two pairs, each with the same number of
endpoints of $\bdd  C$.  (In practice,
$C$ will be the intersection of the boundary of a meridian disk for
$V_2$ with
$P$.)

\begin{lemma} \label{lemma:4punc}

Under these conditions, after a series of edge slides which do not
affect the width of $K \cup \gamma$ or the bridge presentation of $K$,
a meridian for the (possibly different) monotonic edge
$e$  will intersect each component of $C$ in at most one
point.
\end{lemma}

\begin{proof}

Since $e \neq \gamma$, just two of the other $4$ edges of $X$ lie
in $K$.  By edge slides that introduce no new critical points, $e$ can
be shrunk to a point, so that the only vertex in $X$ has one edge
ascending from it and three edges descending.  The mapping class
group of the $4$--punctured sphere is generated by twists around a
pair of simple closed curves each separating a different pair of
punctures and intersecting each other in exactly two points.  In our
context, that means that the entire mapping class group is generated
by twists among the descending triple of edges. Such a twist is
either an edge slide (if one of the strands is in $\gamma$) or just
an isotopy (if both edges are in $K$).

It's easy to see that $\it some$ simple closed curve in $P$ separating
one pair of punctures from another will intersect each component of
$C$ in at most one point. (To see this, note that the condition on
$|\bdd P \cap C|$ guarantees that all arcs of $C$ that have both
ends on a given boundary component of $P$ will be parallel.)  After
an appropriate sequence of braid moves, we can arrange that this
simple closed curve is the meridian of $e$.  Afterwards, unshrink $e$,
recovering a thin presentation of
$K \cup \gamma$ with the required property.  (Note that during this
process the monotonic edge $e$ may switch from being in $K$ to being
in $\gamma$ or vice versa.)
\end{proof}

Of course a symmetric statement holds if there were instead an edge
of $K \cup \gamma$ lying entirely below $S(z)$.

\begin{theorem} \label{theorem:level1}

Let $K$ be a tunnel number one knot in minimal bridge position and
$\gamma$ be a tunnel for $K$. Put $K \cup \gam$ in thin position rel
$K$. Then if
$\gamma$ is an arc, the height function is monotonic on $\gamma$ and
$\gamma$ can be slid and isotoped rel $K$ to be level.  If
$\gamma$ is an eyeglass $\gamma_a \cup \gamma_c$, then the height
function has a single maximum and a single minimum on the circuit
$\gamma_c$, and $\gamma$ can be slid and isotoped, rel $K$, so that
$\gam_c$ is level.
\end{theorem}

\begin{proof}  Put $K \cup \gamma$ in thin position (following the
procedure of Lemma \ref{lemma:4punc} if applicable) so that there is a
level $z$ as in Lemma \ref{lemma:level}.  That is,
$S(z)$ simultaneously has an upper disk $D_u$ and a lower disk
$D_l$, and these disks are disjoint. Moreover, if there is an
edge $e \neq \gamma$ of the graph $K \cup \gamma$ lying
entirely above  (below) $S(z)$, then $\bdd D_u$ (respectively $\bdd
D_l$) intersects the meridian of the edge at most once.

If $\gamma$ itself is disjoint from
$S(z)$ then, because $K \cup \gamma$ is in bridge position, it is
easy to slide
$\gamma$ into $S(z)$, completing the proof.  So we will henceforth
assume that if an edge $e$ of $K \cup \gamma$ is disjoint
from $S(z)$, then $e \neq \gamma$ so, following the application
of Lemma
\ref{lemma:4punc}, a meridian of $e$ will intersect an
upper or lower disk in at most one point.

Let $\aaa_u$ and $\aaa_l$ denote the  (disjoint) arcs
$\bdd D_u \cap S(z)$ and $\bdd D_l \cap S(z)$ respectively.  We will
say that $D_u$ (respectively $D_l$) is {\em bad} if $\aaa_u$ (respectively
$\aaa_l$) is a loop.  Otherwise the disks are {\em good}.  Notice
that if an edge $e$ of $K \cup \gamma$ is disjoint
from $S(z)$, then it is necessarily also disjoint from any bad disk,
since the boundary of a bad disk crosses each meridian an even number
of times.

\begin{definition}

An {\em upper cap} (respectively {\em lower cap}) to $K \cup \gamma$ at
$S(z)$ is an imbedded disk $C \subset S^3$ transverse to $S(z)$ so
that

\begin{enumerate}

\item $C \cap (K \cup \gamma) = \emptyset$,

\item $C \cap S(z) = \bdd C$,

\item $\bdd C$ is essential in $S(z) - (K \cup \gamma)$, and

\item the interior of $C$ lies entirely above (respectively below)
$S(z)$.

\end{enumerate}

\end{definition}

\begin{claim} \label{claim:not2caps}
There cannot be both an upper cap and a lower cap whose
boundaries are disjoint.
\end{claim}

\begin{proof} The boundary of a cap divides $S(z)$ into two disks.  If
the cap $C$ is an upper (respectively lower) cap, then, for each of the
disks in $S(z) -
\bdd C$, an arc of $K - S(z)$ or $\gamma - S(z)$ lying above (respectively
below) $S(z)$ can be isotoped to lie in the disk, indeed just below
(respectively above) $S(z)$.  For upper cap $C_u$ and lower cap $C_l$
with disjoint boundaries, there are disk components of $S(z) -
C_u$ and $S(z) - C_l$ that are disjoint.  Call these disks $S_u$ and
$S_l$.  Then some arc of $K - S(z)$ or $\gamma - S(z)$ lying above
(respectively below) $S(z)$ can be isotoped to lie in $S_u$ (respectively $S_l$).
If in both cases the arc were in $K$ this would violate thin position
of $K$, so at least one of the arcs (say the one in $S_u$) must be in
$\gamma$.  But even then the isotopy would violate thin position of
$K \cup \gamma$ rel $K$, since pushing both arcs through $S(z)$
would reduce $W(K\cup\gamma)$ but would not change the bridge
position of $K$ since $S(z + \epsilon)$ would remain a middle sphere
for $K$.
\end{proof}

\begin{claim} \label{claim:killinterior}
If there is an upper disk (respectively lower disk) and a disjoint lower cap
(respectively upper cap) then we can find such a pair for which the interior
of the upper (respectively lower) disk is disjoint from $S(z)$.
\end{claim}

\begin{proof} Let $B_u$ and $B_l$ denote the balls which are the
closures of the region above $S(z)$ and below $S(z)$ respectively.
It is easy to choose $D_u$ so that each component of $D_u \cap B_u$ is
incompressible in the complement of $K \cup \gamma$.

Since $K \cup \gamma$ is in bridge position, there is a disjoint
collection of disks $\Ddd$ (called {\em descending disks}) so that

\begin{enumerate}

\item $\inter(\Ddd)$ lies above $S(z)$,

\item each boundary component of $\Ddd$ consists of an arc in
$S(z)$ and an arc in $K \cup \gamma$, and

\item the complement of an open neighborhood of $\Ddd \cup K
\cup \gamma$ in $B_u$ is also a $3$--ball.

\end{enumerate}

For example, if no vertex of $K \cup \gamma$ lies above $S(z)$, then
$\Ddd$ is just a set of disks defining a parallelism between the set
of arcs lying above $S(z)$ and arcs on $S(z)$. It is possible to
choose $\Ddd$ so that all intersection points of $\bdd \Ddd$ with
$\bdd D_u$ lie on $S(z)$. This is obvious if the component of $(K
\cup \gamma) - S(z)$ incident to $D_u$ has three ends, since in a
$3$--punctured sphere, arcs connecting two of the punctures to the
third can be isotoped off an arc with both ends on the third.  If the
component $X$ of $(K \cup \gamma) - S(z)$ incident to $D_u$ has four
ends then we have carefully contrived, via Lemma \ref{lemma:4punc},
that $\bdd D_u$ runs at most once along the edge of $K \cup \gamma$
that lies in $X$.  This observation allows us to revert to the case
of the $3$--punctured sphere as above.

The proof of the claim will be by induction on $|D_u \cap
\Ddd|$.  If $|D_u \cap \Ddd| = 0$ then each component of $D_u
\cap B_u$ lies in a ball disjoint from $K \cup \gamma$ so, by
incompressibility, each component is a disk.  Since a
neighborhood of $\bdd D_u$ lies in $B_u$, it follows that
$D_u$ must lie entirely inside $B_u$ as required.

So suppose $D_u \cap \Ddd \neq \emptyset$.  A simple innermost disk
argument could eliminate a closed curve of intersection, so we can
take all components of intersection to be arcs. Surprisingly, we may
also assume that the lower cap is a slight push-off of a disk
component of $D_u \cap B_l$.  Indeed, consider an innermost disk of
$D_u - S(z)$.  If it lies in
$B_u$ then it is an upper cap disjoint from the lower cap which we
know by Claim \ref{claim:not2caps} to be impossible.  If it lies in
$B_l$ then we may as well take a slight push-off as our lower cap.

This surprising fact means that an outermost arc of $D_u \cap
\Ddd$ in $\Ddd$ can be used to $\bdd$--compress $D_u \cap B_u$ to an
arc that is disjoint from the lower cap.  This boundary compression
defines an isotopy on the interior of $D_u$ that reduces $|D_u
\cap \Ddd|$ without disturbing the disjoint lower cap.
Further, this isotopy preserves the incompressibility of
$D_u\cap B_u$ in the complement of $K\cap\gamma$.
After the isotopy, the result follows by induction.
\end{proof}

\begin{claim} \label{claim:capdisk}
There cannot be both an upper disk and a disjoint lower
cap (or a lower disk and a disjoint upper cap).
\end{claim}

\begin{proof}  Following Claim \ref{claim:killinterior} we can assume
that the upper disk $D_u$ lies entirely above $S(z)$.  The boundary
of the lower cap divides $S(z)$ into two disks; let $S_l$ be the one
disjoint from $D_u$.

If $D_u$ is good, then $D_u$ can be used to move
an arc component of  $K - S(z)$ or $\gamma - S(z)$ lying above $S(z)$
down to just below $\aaa_u$. We have already noted that the lower cap
ensures that an arc component of  $K - S(z)$ or $\gamma - S(z)$ lying
below $S(z)$ can be moved up just above $S_l$.  This violates thin
position just as in Claim \ref{claim:not2caps}.

If $D_u$ is bad then the loop $\aaa_u$ divides $S(z)$ into two disks,
let $S_u$ be the one that is disjoint from $S_l$.  Now $\bdd D_u$ is
incident to exactly one edge of $(K \cup \gamma) - S(z)$:  this is
obvious if each component of $(K \cup \gamma) - S(z)$ has at most three
ends and has been previously arranged if a component of $(K \cup
\gamma) - S(z)$ has four ends.  The union of $D_u$ and a neighborhood
of that edge form a disk much like an upper cap. Since $S_u$ contains
some points of  $S(z)\cap (K\cup\gamma)$,  some
edge of $K - S(z)$ or $\gamma - S(z)$ can be slid or isotoped to lie
in $S_u$.  Again this violates thin
position as in Claim \ref{claim:not2caps}.
\end{proof}

\begin{claim} \label{claim:badgood}
There cannot simultaneously be an upper bad disk and a lower good
disk that are disjoint (or, symmetrically, a lower bad disk and an
upper good disk that are disjoint).
\end{claim}

\begin{proof}
We may assume that the interiors of both disks are disjoint from
$S(z)$, for, if not, then an innermost curve of intersection in the
disk would cut off a cap that would contradict Claim
\ref{claim:capdisk}.   Let $D_u$ be the upper bad disk and $D_l$ be
the lower good disk. The loop
$\alpha_u$ divides $S(z)$ into two disks $S_l$ and $S_u$ with
$\alpha_l \subset S_l$. The disk $D_l$ can be used to move
an arc component of  $K - S(z)$ or $\gamma - S(z)$ lying below $S(z)$
up to just below $\aaa_l$.  Since $S_u$ contains some points
of  $S(z)\cap (K\cup\gamma)$, some
edge of $K - S(z)$ or $\gamma - S(z)$ can be slid or isotoped to lie
in $S_u$.  Again this violates thin
position as in Claim \ref{claim:not2caps}.
\end{proof}

By Lemma \ref{lemma:level}, there are an upper disk $D_u$ and a lower
disk $D_l$ for $S(z)$ such that $D_u\cap D_l=\emptyset$.  We may
assume that the interiors of both disks are disjoint from
$S(z)$, for, if not, then an innermost curve of intersection in the
disk would cut off a cap that would contradict Claim
\ref{claim:capdisk}.  By Claim
\ref{claim:badgood} we can assume that either both are good or both
are bad.

If both are good, then they can be used to slide and isotope an arc
component of $K - S(z)$ or $\gamma - S(z)$ lying below $S(z)$ to
just above
$\aaa_l \subset S(z)$ and an arc component of  $K - S(z)$ or $\gamma -
S(z)$ lying above $S(z)$ to be moved just below $\aaa_u$.  If
$\aaa_l$ and $\aaa_u$ have no endpoints in common, then these moves
could be done simultaneously and $K \cup \gamma$ would not have been
thin. The same contradiction arises if they have one end point in
common, unless it is in $\gamma$ and the other two end points are in
$K$.  In that case the simultaneous isotopies and slide make $\gamma$
level as required.  If  $\aaa_l$ and $\aaa_u$ have both end
points in common, they can't both be in $K$ (since $K$ is knotted)
and if both end points are in
$\gamma$ then there is a circuit $\gamma_c \subset \gamma$ that can be
made level, as required. If one common end is in $K$ and the other
in $\gamma$ then $\gamma$ can be moved to a level loop and again we
are done.

If both the upper and lower disks are bad, then the argument is
much the same.  That is, there are disjoint
disks in $S(z)$ bounded by the loops $\aaa_u$ and $\aaa_l$, and arc
components of $K - S(z)$ or
$\gamma - S(z)$ can be moved into those disks.   If these loops were
based at different points, or at a common point lying in $K$, this
move would violate thin position.  If the loops are based at a common
point and it is in $\gamma$, then $\gamma$ can be made level, as
required.
\end{proof}

\begin{definition}

Suppose $K$ is a tunnel number one knot and $\gamma$ is a tunnel for
$K$ in the form $\gamma = \gam_a \cup \gam_c$, where $\gam_c$ is a
circuit.  Suppose $K$ is in minimal bridge position with respect to a
height function $h$,
$\gam_c$ is level with respect to $h$, and $K \cup \gam_a$ is in
Morse position with respect to $h$ with no critical point at the
same height as $\gam_c$.  Then $K \cup \gamma$ is in {\em special
position} with respect to $h$.

\end{definition}

It follows from Theorem \ref{theorem:level1} that for $K$ in thin
position, either the entire tunnel $\gamma$ can be made level, or $K
\cup \gamma$ can be put in special position without affecting the
bridge position of $K$.

We now mimic much of the previous argument:

\begin{definition}

Suppose $K \cup \gamma$ is in special position
(so in particular, $K$ is in minimal bridge position). Let $t_0 <
\ldots < t_n$ be the list of critical heights of $K \cup \gamma_a$,
plus the level of the cycle $\gam_c$.  Let $s_i, 1 \leq i \leq n$
be generic levels chosen so that $t_{i-1} < s_i < t_i$ and let $t_j$
be the level of the vertex $K \cap \gamma_a$ (ie, the vertex {\em
not} on
$\gamma_c$).  Define the width
$W(K \cup \gamma)$ to be the integer $$2(\Sigma_{i \neq j}|S(s_i)\cap
(K
\cup \gamma_a)|) +  |S(s_j)\cap (K \cup \gamma_a)|.$$

For $K$ in minimal bridge position {\it special thin position} of a
pair $(K, \gamma)$ is a special position presentation  which minimizes
$W(K \cup \gamma)$ rel $K$.

Much as in Section \ref{section:thinning}, maxima can be moved past
maxima and minima past minima without affecting the width.  Moving a
maximum below a minimum will always reduce width. What is new is
this:  moving a maximum below the level of $\gamma_c$ or moving a
minimum above the level of $\gamma_c$ will always thin the special
presentation.  This is easily checked.  In particular, the
overarching principle (Proposition \ref{prop:overarch}) remains
intact, and applies also to the level just above $\gamma_c$ (and,
symmetrically, to the level just below $\gamma_c$).

A {\em special bridge position} of $K \cup \gamma$ is a special
position in which every minimum of $K \cup
\inter(\gamma_a)$ occurs at a lower level than every maximum.  (Note
that nothing is said about the level of $\gamma_c$.)
\end{definition}

\begin{proposition}\label{prop:spec.thin=bridge} Suppose that
$(K,\gamma)$ is in special thin position. Then, $K\cup\gamma$ is in
special bridge position.
\end{proposition}

\begin{proof}  The proof is essentially identical to that of
Proposition \ref{prop:thin=bridge}, Case B.  Only easier: the cases in
which $P'$ has a boundary component on $\gamma_c$ don't arise.
\end{proof}

\begin{lemma}  \label{lemma:thin.circuit} Suppose that
$(K,\gamma)$ is in special thin position, so in special bridge
position. Let $x$ be the level of the highest local minimum and  $y$
the level of the lowest local maximum of $K\cup\gamma$.  Then the
level of $\gamma_c$ lies between $x$ and $y$.
\end{lemma}

\begin{proof}  Suppose, on the contrary, that the level of $\gam_c$
is greater than $y$ and choose a level $z$ just below the level
of $\gam_c$ (and so above $y$).  Repeat the argument of Proposition
\ref{prop:thin=bridge}, Case B, at level $z$.  Either a maximum above
the level of $\gam_c$ can be moved down to level $z$, thinning $z$,
or a subarc below $z$ can be moved up to the level of $z$.  Since
below $z$ the next critical level is at a maximum, this necessarily
thins the presentation.
\end{proof}

\begin{lemma}\label{lemma:spec.level1} Suppose that
$(K,\gamma)$ is in special thin position, $x$ is the level of the
highest local minimum,  $y$ is the level of the lowest local
maximum and $\gamma_c$ is at level
$z$ between $x$ and $y$. Then there is not an upper disk for any
$S(t), x < t < z$ or a lower disk for any  $S(t), z <
t < y$.
\end{lemma}

\begin{proof}  The proof mimics that of Lemma \ref{lemma:level}. As
in that lemma, we can arrange that the levels $x$ and $y$ are not
the levels of the other vertex.  Then suppose, say, that $S(t), z < t
< y$ has a lower disk.  We know that $S(y - \varepsilon)$ has an upper
disk, so for some level between $t$ and $y$ there is a level $t_0$
with both an upper and lower disk. The proof is now essentially
identical to that of Theorem \ref{theorem:level1}, except that the bad
lower disk $D_l$ might intersect $K \cup \gamma$ on a component
which is the union of $\gamma_c$ and a subarc of $\gamma_a$.  But
it's easy to see that the arc $\bdd D_l \cap (K \cup \gamma)$ must
then intersect a meridian of $\gamma_c$ exactly once, and so can be
used to slide $\gamma_c$ up to level $t_0$ (and, of course, disjoint
from the slide or isotopy given, as in the proof of Lemma
\ref{lemma:level}, by $D_u$).  The result would be to thin the
presentation.
\end{proof}

\begin{theorem}\label{theorem:spec.thin=level} Suppose that
$(K,\gamma)$ is in special thin position. Then, $\gamma$ may be
isotoped and slid, rel $K$, until it is level.
\end{theorem}

\begin{proof}
Let $z$ be the level of $\gamma_c$, necessarily lying between the
level of the highest minimum and the level of the lowest maximum.

Thicken a neighborhood of
$\gamma_c$ slightly.  Then $S(z) - \eta(\gamma_c)$ consists of two
disks, $S_1$ and $S_2$, each punctured by $K$ and $\gam_a$.  As above,
let $D$ be a meridian disk for $V_2$, with $\bdd D$ on the boundary
of a regular neighborhood of $K \cup \gamma$.  Isotope $D$ to
minimize the pair, in lexicographic order, $(|\bdd D \cap S(z)|,|D
\cap S(z)|)$.  Suppose $D \cap S(z)$ contains closed components, and
let $C$ be an innermost disk cut off by $S(z)$ from $D$.  $C$ is a
cap, lying either above or below $S(z)$.  If, say, it lies above,
then $\bdd C$ bounds a disk in, say, $S_1$.  Then whatever maxima lie
between $S(z)$ and the cap can be pushed down below $S_1$, thinning
the presentation but not affecting the bridge presentation of $K$.  We
conclude that all components of $D \cap S(z)$ are properly imbedded
arcs.

Let $D_u$ be an outermost disk cut off from $D$ by $S(z)$ and, say,
$D_u$ lies above $S(z)$.  Let $\aaa_u$ be the arc $\bdd D_u \cap
S(z)$.  If the ends of $\aaa_u$ lie at distinct points of $(K \cup
\gam_a) \cap S(z)$ then $D_u$ could be used to isotope or slide a
maximum below level $z$, thinning the presentation.  Similarly, if
$\aaa_u$ is a loop with ends at a point in  $K \cup \gam_a$, it
bounds a disk in, say, $S_1$.  Again, whatever maxima have ends in
that disk can be pushed below $S_1$, thinning the presentation.

The possibility remains that one or both ends of $\aaa_u$ lie on
$\bdd S_i$.  In this case, $\gam_a$ must be incident to
$\eta(\gamma_c)$ from above. Consider the annulus half $A$ of
$\bdd\eta(\gamma_c)$ lying above level $z$, punctured by the end of
$\gamma_a$.  No component of $\bdd D \cap A$
can have both ends at the puncture, for otherwise there would be a
lower disk cut off from $D$ just above $\eta(\gamma_c)$, contradicting
Lemma \ref{lemma:spec.level1}. (See Figure 11.) Hence there is a
spanning arc of
$A$ which is disjoint from $\bdd D \cap A$.  This observation allows
us to treat $\gam_c$ much like just another vertex in the previous
argument.  In particular, if $\aaa_u$ has its other end at a point in
$\gam_a \cap S_i$ then the end of $\gam_a$ cut off by that point
could be isotoped just below $S_i$, thinning the presentation.  If
the other end is at a point in $K \cap S_i$ then $\gam_a$ can be made
to lie in an $S_i$, completing the proof.

Finally, suppose that both ends of $\aaa_u$ lie in $\bdd S_i$.  Then $\bdd
D_u$ must run up $\gam_a$.  For otherwise there is an arc of $\bdd D \cap
A$ which does not span $A$, namely $\bdd D_u \cap A$.  Any such
non-spanning arc, either in $A$ or in the annulus $\bdd\eta(\gamma_c) - A$,
must cut off a disk containing the end of $\gam_a$, since $\bdd D \cap
S(z)$ has been minimized.  So there can be no such arc in the annulus
$\bdd\eta(\gamma_c) - A$ and every such arc in $A$ must have its ends on
the same component of $\bdd A$.  The former fact implies that $\bdd D$
intersects both components of $\bdd A$ in the same number of points and the
latter would imply that $\bdd D$ intersects one component more often than
the other, a contradiction.

Having established that $\bdd D_u$ runs up $\gam_a$ note that, just as in
the case of a bad upper disk, $D_u \cup \gam_a$ defines a disk above $S(z)$
which, like an upper cap, can be used to push maxima down to level $S(z)$.
One of these contains $\gam_a$; when it is pushed down to level $S(z)$ the
entire
eyeglass will be level, as required.
\end{proof}

\begin{figure}[ht!]
\cl{\small
\SetLabels 
\E(0.8*0.1)lower disk for higher plane\\
\B(0.43*0.92)$\gamma_a$\\
\E(0.415*0.36)$\gamma_c$\\
\B(0.22*0.15)$A$\\
\endSetLabels 
\AffixLabels{\includegraphics[width=.9\textwidth]{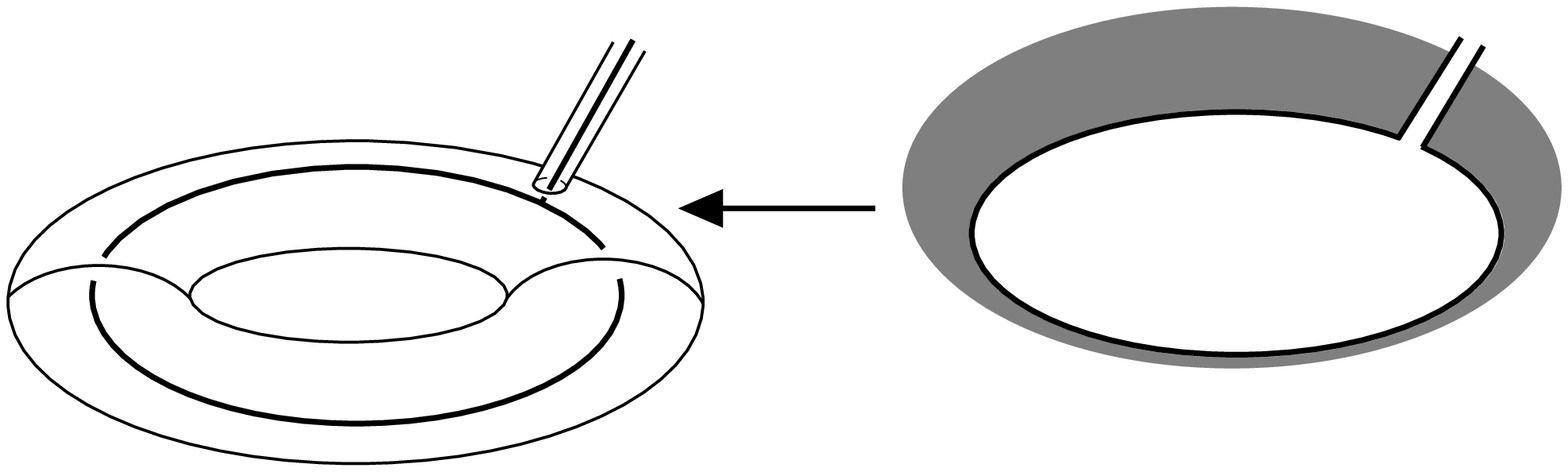}}}
\caption{}
\end{figure}

\section{A characterization of unknotting tunnels}

Let $K$ be an tunnel number one knot of bridge number $n$ and suppose
$\{ B_1,B_2\}$ a pair of 3--balls which gives a minimal bridge
decomposition  of $K$.  That is, $S^{3}=B_{1}\cup B_{2}, B_{1}\cap
B_{2}=\partial B_{1} =\partial B_{2}$ and $(B_i,B_i\cap K)$ is a
trivial tangle of $n$ components
$(i=1,2)$. Let $S=\partial B_1=\partial B_2$ be the middle sphere and
for $i = 1, 2$, let $K_i$ denote the collection of arcs
$B_i\cap K$, parallel to a collection of arcs in $S$.  We have
shown that any tunnel $\gamma$ for $K$ can be slid and isotoped to lie
in $S$.  In this section we glean a bit more information
about how $\gamma$ lies in $S$.  We show:

\begin{theorem} \label{theorem:charact}  For one of
$(B_i,K_i,\gamma)$, $i=1,2$, either:

\begin{enumerate}

\item $\gamma$ is an arc with its ends on different components of
$K_i$ and $K_i$ is parallel to a collection of arcs in $S - \gamma$,

\item  $\gamma$ is an arc with both ends on the
same component of $K_i$.  In this case, $\gamma$ can be slid and
isotoped in $B_i$ so that it lies in $S$ as a loop with
its ends at the same point of $\bdd K_i$, or

\item $\gamma$ is an eyeglass and a disk that $\gamma$ bounds in $S$
contains exactly one end of each of $n-1$ components of $K_i$.

\end{enumerate}

\end{theorem}

\begin{figure}[ht!]
\cl{\small
\SetLabels 
\B(0.17*0.91)$\gamma$\\
\E(0.8*0.6)$\gamma$\\
\T(0.42*0.39)$\gamma$\\
\T(0.22*0.5)(1)\\
\T(0.79*0.5)(2)\\
\T(0.5*-0.015)(3)\\
\endSetLabels 
\AffixLabels{\includegraphics[width=.6\textwidth]{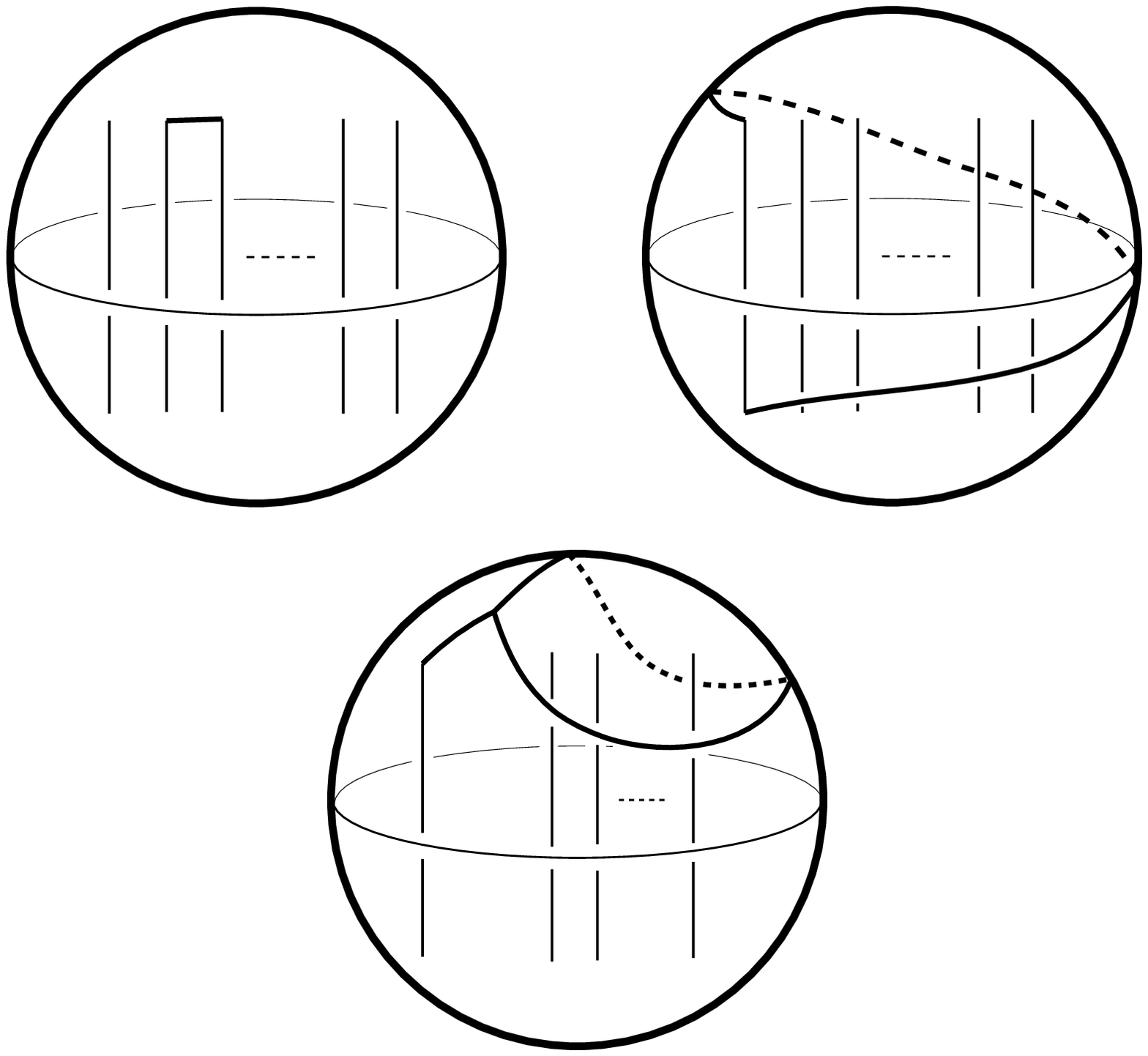}}}
\vglue 2mm\caption{}
\end{figure}

See Figure 12 for typical examples.
In the last case (hence also in the second case), it is shown in
\cite{GOT} that there is a possibly different minimal bridge
position for $K$  with middle sphere $S'$ for which $\gamma$ appears
as part of a `quasi-Hopf tangle' (see Figure 12 (3) or
\cite[Figure 2]{GOT}). In words,
$S'\supset\gamma$ bounds a ball
$B'$ so that the collection of arcs $K \cap B'$ is parallel
to a collection of arcs $\kappa \subset S'$ and the interior of
exactly $n-1$ of the arcs in $\kappa$ intersects $\gamma$, each in
precisely one point.

Figure 13 shows how to construct many examples where it is necessary
to switch minimal bridge positions of $K$ in order to see the
quasi-Hopf tangle.

Thus we have:

\begin{corollary}
By sliding and isotopy of $\gam$, we can find a minimal
bridge sphere of $K$ that bounds a ball in which $K \cup \gam$ is as
seen in Figure 12 (1) or (3).
\end{corollary}

Besides Figure 3, here is another way in which the difference between
the corollary and the theorem can be illustrated.  Imagine replacing
the right-most two vertical arcs in Figure 12 (3) by a rational
tangle, chosen to retain the property that each of the two arc has one end
on each side of the loop formed by $\gam$.  The resulting picture would satisfy
the conclusions of the theorem, but not of the corollary.

\begin{figure}[ht!]
\cl{\small
\SetLabels 
\E(0.42*0.825)full twists\\
\E(0.42*0.295)full twists\\
\E(0.42*0.57)isotopic\\
\E(0.1*0.62)odd \# of half-twists\\
\endSetLabels 
\AffixLabels{\includegraphics[width=.9\textwidth]{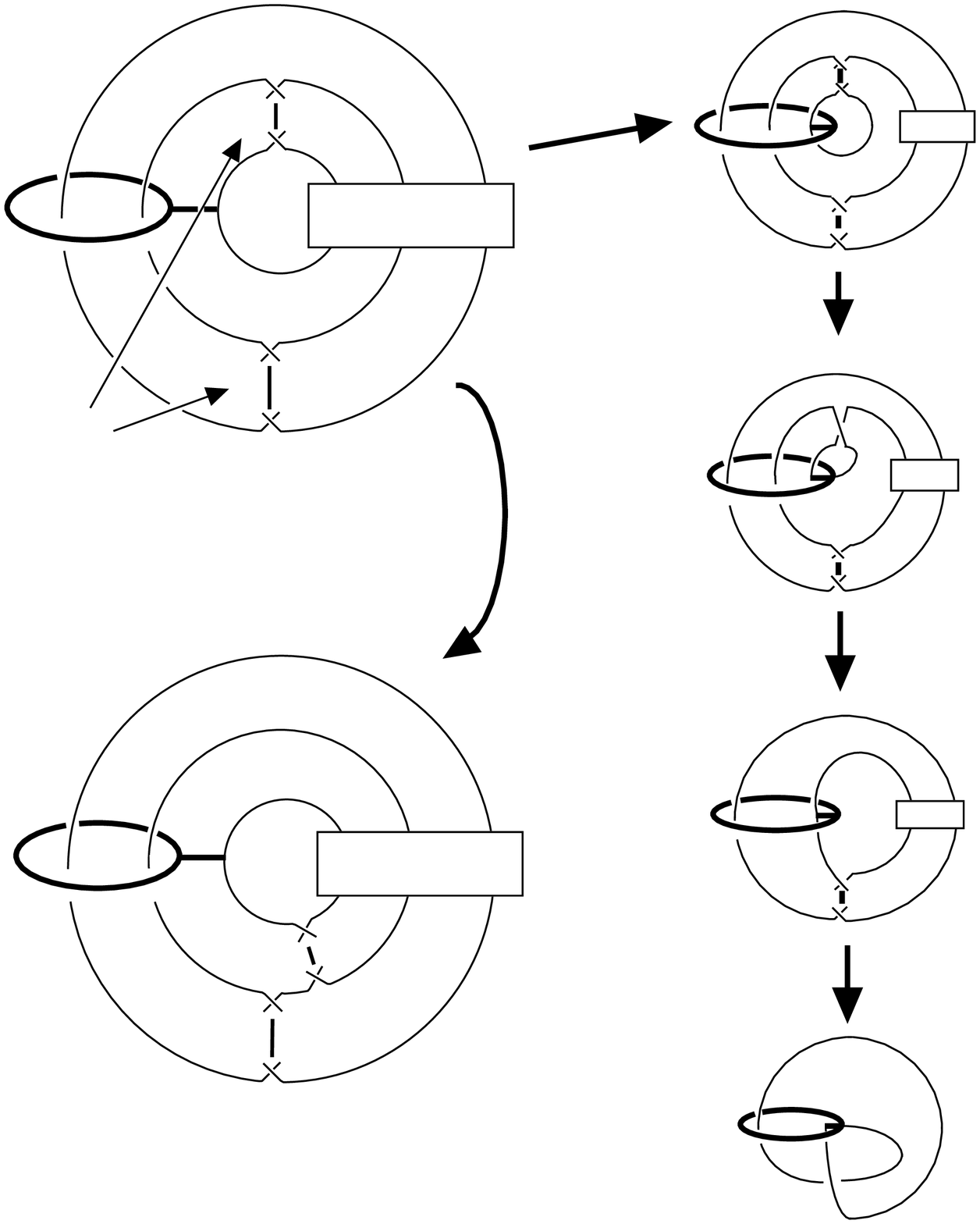}}}
\vglue 2mm\caption{}
\end{figure}

\begin{proof}  The initial work is done above. Recall the
definitions of $V_1, V_2$ and their common boundary $F$, a genus
two surface.  With
$K$ in minimal bridge position we proceed to put $K \cup \gam$ in thin
position rel $K$ and then put $\gam$ into a level sphere that lies
above all the local minima for $K$ and below all the local maxima.
The end(s) of $\gam$ are incident to one or two maximal arcs of $K$
lying above the level of $\gam$ and one or two minimal arcs lying
below.  We call these arcs of $K$ {\it contiguous} to
$\gam$. We divide into three cases.
\medskip

{\bf Case 1}\qua  Some maximal arc (minimal arc) of $K$ not contiguous to
$\gam$ can be pushed below (above) the level of $\gam$ without
changing the bridge presentation of $K$.
\medskip

Push as many maximal arcs as possible below the level of $\gam$.  Let
$S_-$ be a level sphere lying just below $\gam$ but above some maxima
of
$K$. Maximally compress the planar surface
$S_- - \eta(K)$ to get a planar surface $P$ which is
incompressible in
$V_2 = S^3 - \eta(K \cup \gamma)$.  Then, as in Section
\ref{section:thin=bridge}, $P$ is a collection of boundary
parallel planar surfaces in $V_2$ by Proposition \ref{prop:morimoto}.
Much as in Section \ref{section:thin=bridge}  we observe this
contradiction: if an {\em arc} component of $(K \cup \gamma) - S_-$ is
parallel to an arc of $S_-$ this either violates thin position (if
the arc lies below $S_-$) or allows us to push a maximal arc below
$\gam$ (if the arc lies above $S_-$). The latter move might alter the
bridge presentation of $K$, but can be refined using the techniques
of Claim \ref{claim:killinterior} so that it does not, which is a
contradiction to assumption.  Here is a sketch of this refinement:
if there is an upper cap for
$S_-$ then an arc can be pushed down without altering the bridge
position of
$K$.  Since there can be no upper cap, any upper disk either has
interior disjoint from
$S_-$ or is disjoint from some lower cap.  But then it follows from
the argument of Claim \ref{claim:killinterior} that an upper disk can
be found that is disjoint from $S_-$.  So the isotopy it describes
cannot alter the bridge presentation of $K$.

Thus we conclude that
there is no such arc component of $(K \cup \gamma) - S_-$, so
the unique component $\gam_+$ of $(K \cup \gamma) - S_-$ lying above
$S_-$ is the one containing $\gam$.  Moreover, the part $F \cap
\eta(\gam_+)$ of the splitting surface adjacent to $\gam_+$  must be
a planar surface parallel to a component of $P$.  In particular
$\gam_+$ cannot contain a circuit and must be parallel to a subgraph
of $S_-$. So $\gam_+$ is a tree with four ends and has a single
interior edge, necessarily
$\gam$.  This means that $\gam$ is a level edge connecting the
highest two maxima of $K$, matching the first
possibility described in the theorem.
\medskip

{\bf Case 2}\qua  The maximum (minimum) of a contiguous maximal (respectively
minimal) arc of $K$ can be pushed to the level of $\gam$ without
altering the bridge presentation of $K$.
\medskip

Suppose $\gam$ is an arc with ends on two different maximal arcs.
Once one of these arcs is pushed down below $S$, the descending
disks from all the other maximal arcs can be made disjoint from
$\gam$, just by sliding any intersection off the end of $\gam$.  In
particular, a maximal arc not contiguous to $\gam$ (if one exists) can
be pushed below $S$, so we conclude as in Case 1. If there are no
maximal arcs other than the two contiguous to $\gam$, that is in the
$2$--bridge case, the argument still allows us to push both these
maximal arcs to the level of $\gam$, matching the first possibility
described in the theorem.

Suppose $\gam$ is an arc with ends on the same maximal arc $k_0$.
Since the assumption is that $k_0$ can be pushed to the same level as
$\gam$, there is a disk $D_0$ whose boundary consists of the union of
$k_0$ and a subarc of $S - \gamma$.  Slide one end of $\gam$ over
$k_0$ so that $\gam$ becomes a loop, then use $D_0$ to push the new
subarc of $\gam$ down onto $S$. This move demonstrates that $\gam$
satisfies the second possibility described in the theorem.

Finally, suppose that $\gam$ is an eyeglass, which we may as well
take to be a loop based at a point of $K \cap S$ (ie, set
$\gam_a = \emptyset$).  Once the maximum of its contiguous maximal arc
is pushed down to the level of $\gam$, we can assume that the
interior of one of the two disks bounded by $\gam$ intersects $K$ in
$k \leq n-1$ points, where $n$ is the bridge number of $K$, ie, the
number of arcs in $K_i$.  It then follows from \cite{GOT} that a
regular neighborhood of that closed disk has boundary a middle sphere
$S'$ for a possibly different bridge presentation of $K$.  But $|S'
\cap K| = 2k + 2$ and since $n$ is the minimal bridge number of $K$
we have $2k + 2 \geq 2n$.  Hence in fact $k = n-1$, matching
description 3 of the theorem.
\medskip

{\bf Case 3}\qua  The general case.
\medskip

$V_1\cap S$ consists of some meridian
disks $E_i$ of $K$ together with either a disk
component $E_0$ or an annulus component $A_0$ that cuts
$\gam$ lengthwise.  It is a disk when $\gam$ is an arc
and an annulus when $\gam$ is an eyeglass.

Let $D$ be a meridian disk for $V_2$ and $P$ be the planar surface $S
- V_1$.  If $D \cap P$ has a closed component, consider
an innermost disk component of $D - P$; if not, consider an
outermost disk of $D - P$, cut off by an outermost arc.  The argument
is similar in both cases, so we focus on a disk $D_0$, an
outermost disk of $D - P$.  With no loss of generality, assume
$D_0$ lies entirely above $S$.
If the arc $\bbb
= F \cap \bdd D_0$ lies on an annulus component of $F - S$ then $D_0$
provides an isotopy of a maximal arc of $K$ to the level
of $S$, hence to a level below $S$ and we conclude as in Case 1.
Otherwise $\bbb$ lies on a component $F_{\gam}$ of $F - S$
neighboring $\gam$.  There are three possibilities:
\medskip

{\bf Subcase 3a}\qua  $F_{\gam}$ is a once-punctured torus, with
$\bdd F_{\gam} = \bdd E_0$.
\medskip

This case arises when $\gam$ is an arc with ends coincident to the
same maximal arc $k_0$ of $K$. The
torus $F_{\gam} \cup (S - E_0)$ bounds a solid torus, namely the union
of a neighborhood of $k_0$ with the ball below $S$.  $\bdd D_0$ is an
essential circle on the boundary of this torus bounding a disk in its
complement.  Hence it is a longitude of the solid torus.  It follows
that $\bbb$ runs exactly once over the arc $k_0$, so in fact $D_0$
provides an isotopy of $k_0$ to a subarc of $S$ disjoint from
$\gam$. Then we conclude as in Case 2.

\medskip
{\bf Subcase 3b}\qua  $F_{\gam}$ is a three-punctured sphere, whose
boundary consists of  $\bdd E_0$ and the boundary of two meridia of
$K$.
\medskip

Suppose first that both ends of $\bbb$ lie on $\bdd E_0$.  Then one
way to think of $D_0$ is as a disk disjoint from $K$, lying above $S$,
and crossing $\gam$ exactly once.  But this implies that the
descending disks $\Ddd$ from the maxima of $K$ can be made disjoint
from $\gam$:  any intersection with $\gam$ can be piped to $D_0$ and
thereby removed.  In particular, all maximal arcs of $K$ can be
isotoped to lie on $S - \gam$ so we conclude as in Case 2.

If exactly one end of $\bbb$ lies on $\bdd E_0$, then $D_0$ can be
used to push a maximal arc contiguous to $\gam$ to the
level of $\gam$ and again we conclude as in Case 2.

Finally, if neither end of $\bbb$ lies on $\bdd E_0$ then a subdisk of
$D_0$ is a lower disk for a level sphere $S_+$ lying a bit above
$S$.  So at some level $S(z)$ above $S$ and below a contiguous maximum
there is both an upper disk $D_u$ and a lower disk $D_l$ cut off by
$S(z)$; these disks are parallel to  subdisks of $D_0$.  Normally, the
existence of an upper and lower disk would not be particularly useful
information, since the component $\gam_+$ of $K \cup \gamma$ lying
below $S_+$ and containing $\gamma$ has four ends, so the arc of
$D_l$ incident to $K \cup \gamma$ might be a complicated arc on a
neighborhood of $\gam_+$.  But here we know that
$\partial\eta(\gam_{+})\cap\partial D_l\subset\beta$,
so $\bdd D_l$ is disjoint from $\bdd E_0$.
We note that $D_u$ is good.

In particular, if $D_l$ is good, then $\bbb$ runs once
across a meridian of $\gam$ and can be used to isotope
$\gam$ up at the same time that we use $D_u$ to bring a
contiguous maximum down. Once this is done, we conclude as in Case 2.
Suppose $D_l$ is bad so the loop
$\bdd D_l$ divides $S(z)$ into two disk components $S_l$ and
$S_u$, with $(\bdd D_u \cap S(z)) \subset S_u$.
Then it is not hard to use
$D_l$ to isotope $\gam$ up to a subarc of $S_l$ and use $D_u$ to
isotope a maximum down to a subarc of $S_u$. Again we conclude as in
Case 1 or 2.
\medskip

{\bf Subcase 3c}\qua  $F_{\gam}$ is a three-punctured sphere, whose
boundary consists of
$\bdd A_0$ and the boundary of  a single meridian of $K$.
\medskip

Consider an outermost disk $D_0$ as above, with $\bbb = \bdd D_0 \cap
F_{\gam}$.  Its ends can't lie on different components of $\bdd A_0$
since these bound disjoint disk components of $S - A_0$.  If one end
lies on a component of $\bdd A_0$ and the other on a meridian of $K$,
then $D_0$ defines an isotopy of a contiguous maximum down to $S$,
and we conclude as in Case 2. If both ends of $\bbb$ lie on a meridian
of $K$, then $\bbb$ must be an arc that crosses a meridian of $\gam$
exactly once.  This situation was dealt with in Lemma
\ref{lemma:spec.level1} above:  A level sphere just above $S$ cuts
off a lower disk from $D_0$, so at some level there will be an upper
disk and a lower disk that allow
$\gam$ and its contiguous maximal arc to be pushed to the same level.
Again, the conclusion follows from Case 2.

Finally, suppose both ends of
$\bbb$ lie on the same component of $\bdd A_0$.  The argument is
much as in that part of Subcase 3b in which both ends of $\bbb$
lie on $\bdd E_0$:  One can think of
$D_0$ as a disk disjoint from
$K$, lying above $S$, and crossing $\gam$ in exactly one point,
a point in $\gam_a$. (Here we drop the conceptual convention that
$\gam_a =
\emptyset$.)  Then $\bdd D_0$ divides $S$ into two disks; let $S_0$
be the one that does not contain $\gam_c$.  Then, as in Subcase 3b,
all descending disks of maxima with ends in $S_0$ can be
made disjoint from $\gam$, first by removing all intersections of
the descending disks with
$D_0$ (thereby ensuring that the descending disks never intersect
$\gam_c$) then by using $D_0$ to pipe away any intersections with
$\gam_a$.  Once the descending disks are disjoint from $\gam$, the
descending disk for the maximum of $K$ contiguous to $\gam$ can be
used to push the maximum to $S$, and the conclusion again follows
from Case 2.
\end{proof}

\section{A brief remark on tunnel number $1$ links.}

The central results here can easily be extended to tunnel number one
links, though the trajectory of proof is a bit different.  It is not
true that the complement of a tunnel number link has no
incompressible meridional planar surfaces.  But if there is such a surface,
then Gordon and Reid show there is one that cuts off a {\em Hopf
tangle}  (see \cite[Figure 12]{GR}) and which has no more boundary
components than the original one. Unknotting tunnels for such links have
been extensively
studied in \cite{GOT}.

In our context, the results of \cite{GR} and \cite{GOT} imply that if $L$
is a tunnel number one link and a minimal bridge position is not
thin, then one way to put $L$ into  minimal bridge position is via
decomposing it into the union of the $n$--string quasi Hopf tangle and the
$n$--string trivial tangle.  (The minimal bridge number of $L$ is then
$n+1$.) In particular, the link can be described as in (3) of Theorem
\ref{theorem:charact}, except that $\gam_c$ is here the unknotted component
of the link.

If a minimal bridge position for the link $L$ is thin, then we can
proceed (as in the case of a tunnel number one knot) to deduce that not
only is the link $L$ in bridge position, but so is the graph $L \cup
\gamma$; indeed, $\gamma$ must lie in a level sphere.

\rk{Acknowledgements}Part of this work was carried out while the
first author was visiting at University of California, Davis. He would
like to express sincere thanks to the department for its hospitality.
The second author was partially supported by an NSF grant and the
Miller Foundation.  The third author was partially supported by an NSF
grant.

\end{document}